\def\date{19.02.08} 

\documentclass[12pt]{article}

\usepackage{amssymb}
\usepackage{amsmath}
\input{liemacs.sty} 

\begin{document}

\renewcommand{\R}{{\mathbb R}}
\renewcommand{\C}{{\mathbb C}}
\newcommand{\bL}{{\mathbb L}}

\def\NN{\mathbb N}
\def\ZZ{\mathbb Z}
\def\QQ{\mathbb Q}
\def\RR{\mathbb R}
\def\CC{\mathbb C}
\def\HH{\mathbb H}
\def\DD{\mathbb D}
\def\TT{\mathbb T}
\def\SS{\mathbb S}

\def\N{\mathbb N}
\def\Z{\mathbb Z}
\def\Q{\mathbb Q}
\def\R{\mathbb R}
\def\C{\mathbb C}
\def\H{\mathbb H}
\def\D{\mathbb D}
\def\T{\mathbb T}
\def\1{{{\bf 1}}}
\def\0{{{\bf 0}}}

\newcommand\ta{\tau}

\newcommand\resp{resp.\ }
\newcommand\ie{i.e., }
\newcommand\cf{cf.\ }

\newcommand \al{\alpha}
\newcommand\be{\beta}
\newcommand\ga{\gamma}
\newcommand\de{\delta}
\newcommand\ep{\varepsilon}
\newcommand\ze{\zeta}
\newcommand\et{\eta}
\renewcommand\th{\theta}
\newcommand\io{\iota}
\newcommand\ka{\kappa}
\newcommand\rh{\rho}
\newcommand\ro{\rho}
\newcommand\si{\sigma}
\newcommand\ph{\varphi}
\newcommand\ps{\psi}
\newcommand\om{\omega}
\newcommand\Om{\Omega}
\newcommand\Ga{\Gamma}
\newcommand\De{\Delta}
\newcommand\Th{\Theta}
\newcommand\La{\Lambda}
\newcommand\Si{\Sigma}

\newcommand{\proj}{\mathop{{\rm pr}}\nolimits}
\newcommand{\dend}{\mathop{{\rm dend}}\nolimits}

\newcommand\kk{\mathfrak k}
\newcommand\X{\mathfrak X}
\newcommand\x{\times}
\newcommand\LL{\mathfrak L}

\title{An abstract setting for hamiltonian actions} 
\author{Karl-Hermann Neeb and Cornelia Vizman}

\maketitle

\begin{abstract} 
In this paper we develop an abstract setup for hamiltonian group 
actions as follows: Starting with a continuous $2$-cochain 
$\omega$ on a Lie algebra $\fh$ with values in an $\h$-module $V$, 
we associate subalgebras 
$\sp(\h,\omega) \supeq \ham(\h,\omega)$ of symplectic, resp., 
hamiltonian elements. Then $\ham(\h,\omega)$ has a natural central extension 
which in turn is contained in a larger abelian extension of 
$\sp(\h,\omega)$. In this setting, 
we study linear actions of 
a Lie group $G$ on $V$ which are compatible with a homomorphism 
$\g \to \ham(\h,\omega)$, \ie abstract hamiltonian actions, 
corresponding central and abelian extensions 
of $G$ and momentum maps $J \: \g \to V$. 

\nin {\it Keywords}: central extension, momentum map, hamiltonian action, 
abelian extension, infinite dimensional Lie group \\ 
\nin {\it MSC}: 17B56, 35Q53 
\end{abstract}

\section*{Introduction.}

In \cite{Br93} Brylinski describes how to associate to a 
connected, not necessarily finite-dimensional, 
smooth manifold $M$, endowed with a closed $2$-form 
$\omega \in \Omega^2(M,\R)$, 
a central extension of the Lie algebra $\ham(M,\omega)$ of 
hamiltonian vector fields on $M$. 
If there exists an associated pre-quantum bundle $P$ with connection 
$1$-form $\theta$ and curvature $\omega$ (which is the case if 
$\omega$ is integral and $M$ is smoothly paracompact), 
then Kostant's central extension (\cite{Ko70}) is given by the short exact 
sequence 
\begin{equation}\label{kostant}
\1\to \T\into\Aut(P,\th)_0\onto\Ham(M,\om)\to \1
\end{equation}
of groups, where $\Aut(P,\th)_0$ is the group of those connection preserving 
automorphisms of $P$ isotopic to the identity 
and $\Ham(M,\om)$ is the group of hamiltonian diffeomorphisms of $M$. 
In general, neither $\Ham(M,\omega)$ nor $\Aut(P,\theta)_0$ carries a 
Lie group structure if $M$ is not assumed to be compact. 

However, we have shown in \cite{NV03} that central Lie group extensions 
can be obtained as follows, even if $M$ is infinite dimensional. 
Let $Z$ be an abelian Lie group of the form $Z = \z/\Gamma_Z$, where 
$\Gamma_Z$ is a discrete subgroup of the Mackey complete space $\z$, 
$\omega \in \Omega^2(M,\z)$ a closed $2$-form, and 
$(P,\theta)$ a corresponding $Z$-pre-quantum bundle, i.e., $P$ is a 
$Z$-principal bundle and $\theta \in \Omega^1(P,\z)$ a connection $1$-form 
with $\dd \theta = q_P^*\omega$, where $q_P \: P \to M$ 
is the bundle projection. 
We call a smooth action of a (possibly infinite dimensional) 
connected Lie group $G$ on a (possibly infinite dimensional) 
manifold $M$ {\it hamiltonian} if the derived 
homomorphism maps into hamiltonian vector fields: 
$$\zeta \: \g \to \ham(M,\omega) := \{ X \in {\cal V}(M) \: 
(\exists f \in C^\infty(M,\R))\ i_X \omega = \dd f\}. $$
Then the pullback of $\Aut(P,\theta)$ defines a central Lie group extension
\begin{equation}\label{eq:centext} 
\1\to Z\into\hat G_{\rm cen}\onto G\to \1, 
\end{equation}
i.e., $\hat G_{\rm cen}$ carries a Lie group structure for which 
it is a principal $Z$-bundle over~$G$. A Lie algebra 2-cocycle
for the associated central Lie algebra extension $\hat\g$ of $\g$ by $\R$ 
is given by 
$(X,Y)\mapsto -\om(\ze(X),\ze(Y))(m_0)$ for any fixed element $m_0 \in M$.

The main point of the present paper is to provide an abstract setting 
for this kind of hamiltonian group actions, momentum maps, 
and the associated central Lie group actions. As we shall see in the examples 
in Section~\ref{sec:1}, our setting 
is sufficiently general to cover various kinds of examples 
of different nature described below. 

Starting with a continuous $2$-cochain 
$\omega$ on a Lie algebra $\fh$ with values in a topological $\h$-module $V$, 
we associate the subalgebras 
\begin{align*}
&\sp(\h,\omega) 
:= \{ \xi \in \h \: {\cal L}_\xi\omega = 0, i_\xi \dd_\h \omega = 0\} 
= \{ \xi \in \h \: \dd_\h(i_\xi\omega) = i_\xi \dd_\h \omega = 0\} \\
\supeq &\ham(\h,\omega) := \{ \xi \in \h \: i_\xi\dd_\h\omega=0, 
(\exists v \in V)\ i_\xi\omega = \dd_\h v\}. 
\end{align*}
of {\it symplectic}, resp., {\it hamiltonian elements} of $\h$. 
Then $\ham(\h,\omega)$ has a natural central extension 
$$ \hat\ham(\h,\omega) := \{ (v,\xi) \in V \times \ham(\h,\omega) \: 
\dd_\h v = i_\xi \omega\} $$
by the trivial module $V^\h$ 
which in turn is contained in the larger abelian extension of 
$\sp(\h,\omega)$ by $V$ defined by $\omega$. To obtain the Lie bracket 
on $\hat\ham(\h,\omega)$, we observe that the space 
$V_\omega := \{ v\in V \: (\exists \xi \in \ham(\h,\omega)) \ 
\dd_\h v = i_\xi \om\}$ of so-called 
{\it admissible elements} carries a Lie bracket 
analogous to the Poisson bracket: 
$$ \{v_1, v_2\} := \omega(\xi_2, \xi_1) \quad \mbox{ for } \quad 
\dd_\h v_j  = i_{\xi_j} \omega, $$
and $\hat\ham(\h,\omega)$ is a subalgebra of the Lie algebra direct sum 
$V_\omega \oplus \ham(\h,\omega)$. 

The classical example is given by $\h = {\cal V}(M)$, $V = C^\infty(M,\R)$ 
and a closed $2$-form $\omega$. In a similar spirit is the 
example arising from a Poisson manifold 
$(M,\Lambda)$, where we put $\h = \Omega^1(M,\R)$, 
endowed with the natural Lie bracket (cf.\ Example~\ref{ex:1.5}), 
$V = C^\infty(M,\R)$  
and $\omega := \Lambda$, considered as a $V$-valued $2$-cocycle on $\h$. 
Of a different character are the examples obtained with 
$\h = {\cal V}(M)$, 
$V := \oline\Omega^p(M,\R) := \Omega^p(M,\R)/\dd \Omega^{p-1}(M,\R)$ 
and $\omega(X,Y) := [i_Y i_X \tilde\omega]$, where 
$\tilde\omega$ is a differential $(p+2)$-form on $M$. There are more examples 
arising from associative algebras in the spirit of non-commutative geometry.

In Section~\ref{sec:2} we turn to {\it abstract hamiltonian actions} of
a Lie group $G$ on $V$, \ie actions which are compatible with a homomorphism 
$\zeta \: \g \to \ham(\h,\omega)$ and the $\h$-action on $V$. Then the 
pullback extension  
$\hat\g_{\rm cen} := \zeta^*\hat\ham(M,\omega)$ defines a central extension 
of $\g$ by $V^\h$, and a {\it momentum map} is a continuous linear map 
$J \: \g \to V$ satisfying 
$$ \dd_\h(J(X)) = i_{\zeta(X)}\omega \quad \mbox{ for } \quad X \in \g, $$
i.e., $J$ defines a continuous linear 
section $\sigma := (J,\zeta,\id_\g) 
\: \g \to \hat\g_{\rm cen}$ of the associated 
central extension.  The obstruction to the existence of an equivariant 
momentum map is measured by the $1$-cocycle 
$\kappa \: G \to C^1(\g,V^\h), \kappa(g) := g.J - J$
which in turn can be used to describe the adjoint action of the 
group $G$ on the extended Lie algebra $\hat\g_{\rm cen}$. 

In Section~\ref{sec:3} we briefly discuss the existence of a 
central Lie group extension $\hat G_{\rm cen}$ with Lie algebra 
$\hat\g_{\rm cen}$. Once a momentum map $J$ is given for a hamiltonian 
$G$-action, the only condition 
for the existence of such an extension is the discreteness of the image 
$\Pi_\omega$ 
of a period homomorphism $\per_{\omega_\g} \: \pi_2(G) \to V^\fh$, associated
to the Lie algebra cocycle $\omega_\g := \zeta^*\omega$. 
If $\Pi_\omega$ is discrete, 
there also is an abelian Lie group extension 
$\hat G_{\rm ab}$ of $G$ 
integrating the Lie algebra extension $\hat\g_{\rm ab} := 
V\rtimes_{\om_\g}\g$.

In many situations the period map $\per_{\omega_\g}$ 
is quite hard to evaluate, but one may 
nevertheless show that its image is discrete. To illustrate this point, 
we take in Section~\ref{sec:4} a closer look at a smooth action 
of a connected Lie group $G$ on a manifold $M$, leaving a closed 
$\z$-valued $2$-form $\omega$ invariant. If the group 
$S_\omega = \int_{\pi_2(M)} \omega$ of spherical periods of $\omega$ 
is discrete in $V$, then $Z := \z/S_\omega$ is an abelian Lie group, 
and there exists a Lie group extension 
$$ \1 \to Z \to \hat G_{\rm cen} \to \tilde G \to \1 $$
integrating the Lie algebra $\hat\g_{\rm cen}$. Here $q_G \: \tilde G \to G$ 
denotes a simply connected covering group of $G$, so that 
$\hat G_{\rm cen}$ can also be viewed as an extension of $G$ by a 
$2$-step nilpotent Lie group $\hat\pi_1(G)$ which is a central extension 
of the discrete group $\pi_1(G) = \ker q_G$ by $Z$. 
Here we 
use that the universal covering manifold $q_M \: \tilde M \to M$, endowed 
with $\tilde\omega := q_M^*\omega$, permits us 
to embed $\sp(M,\omega)$ into $\sp(\tilde M, \tilde\omega) 
= \ham(\tilde M,\tilde\omega)$, so that we actually obtain an abstract
hamiltonian action of the universal covering group 
$\tilde G$ for the module $V := C^\infty(\tilde M,\z)$ 
and the central extension $\hat G_{\rm cen}$ acts by quantomorphisms on 
a $Z$-principal bundle $P$ over $\tilde M$. This further leads 
to an abelian Lie group extension 
$$ \1 \to C^\infty(M,Z)_0 \to \hat G_{\rm ab} \to \tilde G \to \1, $$
integrating the cocycle $\zeta^*\omega \in Z^2(\g,C^\infty(M,\z))$. 
In the short Section~\ref{sec:5} we formulate the algebraic 
essence of the Noether Theorem in our context and the paper concludes 
with an appendix recalling some of the results from \cite{Ne04} on 
the integrability of abelian Lie algebra extensions to Lie group 
extensions. 


\section{A Lie algebraic hamiltonian setup}  \label{sec:1}

Let $V$ be a topological module of the topological Lie algebra $\h$,  
i.e., the module operation $\h \times V \to V, (\xi,v) \mapsto \xi . v,$ 
is continuous. We write $(C^\bullet(\h,V),\dd_\h)$ for the 
Chevalley--Eilenberg complex of continuous Lie algebra cochains 
with values in $V$, $Z^p(\h,V)$ denotes the space of $p$-cocycles, 
$H^p(\h,V)$ the cohomology space etc. We further write 
${\cal L}_\xi$ for the operators defining the natural action of $\h$ 
on the spaces $C^p(\h,V)$. 
We refer to \cite{FF01} for the basic notation concerning 
continuous Lie algebra cohomology. 

We first recall 
some results from \cite{Ne06a} 
which we specialize here for $2$-cochains.
Fix a continuous $2$-cochain $\om\in C^2(\h,V)$. 
Then 
$$
\sp(\h,\om):=\{\xi\in\h: {\cal L}_\xi\om=0, i_\xi \dd_\h\omega=0\} 
= \{ \xi \in \h \: \dd_\h(i_\xi\omega)=0, i_\xi\dd_\h\omega=0\} 
$$
is a closed subspace of $\h$ and from 
$$ [{\cal L}_\xi, {\cal L}_\eta] = {\cal L}_{[\xi,\eta]} 
\quad \mbox{ and } \quad 
[{\cal L}_\xi, i_\eta] = i_{[\xi,\eta]} $$
it follows that $\sp(\h,\omega)$ is a Lie subalgebra of $\h$, 
called the Lie algebra of {\it symplectic elements} of $\h$. 
Since $\dd_\h \omega$ vanishes on $\sp(\h,\omega)$, the 
restriction of $\omega$ to this subalgebra is a Lie algebra 
$2$-cocycle. 
Using the Cartan formula ${\cal L}_\xi = i_\xi \circ \dd_\h + \dd_\h \circ i_\xi$, 
we find 
for $\xi,\eta\in\sp(\h,\om)$ the relation 
\begin{equation}\label{dh}
i_{[\xi,\eta]}\om
= [{\cal L}_{\xi},i_\eta]\om 
= {\cal L}_{\xi}(i_\eta\om)
= \dd_\h(i_\xi i_\eta\omega) 
=\dd_\h(\om(\eta,\xi)), 
\end{equation}
showing that the {\it flux homomorphism},  
$$ f_\om: \sp(\h,\om)\to H^1(\h,V), \quad 
\xi \mapsto [i_\xi\om], $$
is a homomorphism of Lie algebras if $H^1(\h,V)$ is endowed with the trivial Lie 
bracket. 
Its kernel is the ideal 
$$ \ham(\h,\om):=\{\xi\in\sp(\h,\omega)
:(\exists v \in V)\ i_\xi\om=\dd_\h v\} $$
of {\it hamiltonian elements}. 
The set of {\it admissible vectors} 
$$ V_\om:=
\{v\in V: (\exists \xi \in \ham(\h,\omega))\ \dd_\h v=i_\xi\om\}$$
contains the subspace 
$$
V^\h=\{v\in V: (\forall \xi \in \h)\ \xi. v=0\} = H^0(\h,V) $$ 
of $\h$-invariant vectors. We then have a well-defined linear map 
$$ q:\ham(\h,\om)\to \dd_\h V_\omega \cong V_\om/V^\h, 
\quad q(\xi)=i_\xi\om $$ 
whose kernel is the radical 
$\rad(\omega,\dd_\h\omega) := \{ \xi \in \h \: i_\xi \omega = 0, 
i_\xi \dd_\h\omega 
= 0\}$ of $\omega$ and $\dd_\h\omega$.  

\begin{proposition} \label{prop:1} 
{\rm(i)} The space $V_\om$ of admissible vectors carries a Lie algebra 
structure defined by 
$$
\{v_1,v_2\}:=\xi_1. v_2=-\xi_2. v_1=-\om(\xi_1,\xi_2) 
\quad \mbox{ for } \quad \dd_\h v_j=i_{\xi_j}\om, j=1,2, $$ 
and for which $V^\h$ is central. We have the following 
exact sequence of Lie algebras 
\begin{eqnarray}\label{hamhomega}
\0\to \rad(\omega,\dd_\h\omega)\into \ham(\h,\om)\stackrel{q}{\onto} \dd_\h(V_\om)\to \0.  
\end{eqnarray}

{\rm(ii)} If $\h_\omega :=  \{ \xi \in \h \: i_\xi \dd_\h \omega = 0, 
[\xi,\rad(\omega,\dd_\h\omega)] \subeq \rad(\omega,\dd_\h\omega)\}$ 
is the normalizer of $\rad(\omega,\dd_\h\omega)$, then $\sp(\h,\omega) \subeq \h_\omega$ 
and the set 
$$ C^1(\h,V)_\omega := \{ i_\xi \omega \: \xi \in \h_\omega\} \subeq C^1(\h,V) $$
inherits a  Lie algebra structure, defined by 
$$ [i_{\xi_1} \omega, i_{\xi_2} \omega] := i_{[\xi_1,\xi_2]}\omega.  $$
For $\alpha_j = i_{\xi_j} \omega$, $j=1,2$, we then have 
\begin{eqnarray}\label{alpha} 
[\alpha_1,\alpha_2] = {\cal L}_{\xi_1}\alpha_2 - i_{\xi_2}{\cal L}_{\xi_1}\omega. 
\end{eqnarray}
The maps 
$$ i^\omega: \h_\omega \to C^1(\h,V)_\omega, \ \  \xi \mapsto i_\xi \omega 
\quad \mbox{ and } \quad  \dd_\h \: V_\omega \to C^1(\h,V)_\omega,\ \  v \mapsto \dd_\h v $$
are homomorphisms of Lie algebras. In particular
\begin{eqnarray}\label{homega}
\0\into \rad(\omega,\dd_\h\omega)\to \h_\om\stackrel{i^\om}{\onto} C^1(\h,V)_\om\to \0
\end{eqnarray}
is an exact sequence of Lie algebras extending~\eqref{hamhomega}. 
\end{proposition}

\begin{proof} (i) First we note that $\{v_1, v_2\}$ is well-defined 
because the formula 
$$ \xi_1.v_2 = i_{\xi_1}\dd_\h v_2 = i_{\xi_1}i_{\xi_2}\omega 
= \omega(\xi_2, \xi_1) $$
shows that each choice of $\xi_2$ with $i_{\xi_2}\omega = \dd_\h v_2$ 
leads to the same value of $\omega(\xi_2, \xi_1)$, and a similar argument 
applies to $\xi_1$. 

By \eqref{dh}, $\dd_\h\{v_1,v_2\}=i_{[\xi_1,\xi_2]}\om$, so that 
$V_\omega$ is closed under the bracket $\{\cdot,\cdot\}$ and 
$q$ is compatible with the brackets. It remains to verify the Jacobi 
identity in $V_\omega$:  
For $\xi_j\in\ham(\h,\omega)$ 
with $\dd_\h v_j=i_{\xi_j}\om$, $j=1,2,3$, we  have 
$\xi_j \in \sp(\fh,\omega)$, so that 
\begin{align*}
0&=(\dd_\h\om)(\xi_1,\xi_2,\xi_3)=\sum_{\rm cycl.}\xi_1.\om(\xi_2,\xi_3)
-\sum_{\rm cycl.}\om([\xi_1,\xi_2],\xi_3)\\
&\stackrel{\eqref{dh}}{=}-2\sum_{\rm cycl.}\om([\xi_1,\xi_2],\xi_3)
=2\sum_{\rm cycl.}\{\{v_1,v_2\},v_3\}.  
\end{align*}

(ii) The inclusion $\sp(\h,\omega) \subeq \h_\omega$ follows from 
\eqref{dh}, so that $\rad(\omega,\dd_\h\omega) \subeq \h_\omega$ is an ideal of 
$\h_\omega$. Hence the set 
$C^1(\h,V)_\omega \cong \h_\omega/\rad(\omega,\dd_\h\omega)$ 
inherits a quotient Lie algebra 
structure for which $i^\omega$ is a morphism of Lie algebras. 
Further, $\dd_\h \: V_\omega \to C^1(\h,V)_\omega$ is a homomorphism 
of Lie algebras
since for $\dd_\h v_j = i_{\xi_j} \omega$, $j=1,2$, equation 
\eqref{dh} leads to 
$\dd_\h\{v_1,v_2\}= i_{[\xi_1, \xi_2]}\omega = [i_{\xi_1}\om,i_{\xi_2}\om]
=[\dd_\h v_1,\dd_\h v_2].$
Now \eqref{alpha} follows from 
$[\alpha_1, \alpha_2] = i_{[\xi_1, \xi_2]}\omega 
= [{\cal L}_{\xi_1}, i_{\xi_2}]\omega  
= {\cal L}_{\xi_1}(i_{\xi_2}\omega) - i_{\xi_2} {\cal L}_{\xi_1}\omega.$
\end{proof}

The pullback by $q$ of the central extension 
$V^\h\into V_\om\onto\dd_\h V_\om$
provides a central extension  
\begin{equation}\label{hatham}
\0\to V^\h\into\widehat\ham(\h,\om) \onto\ham(\h,\om)\to \0, 
\end{equation}
where 
$$
\widehat\ham(\h,\om)=\{(v,\xi)\in V_\om\x\ham(\h,\omega) :i_\xi\om=\dd_\h v\} $$
is endowed with the Lie bracket 
$$[(v_1,\xi_1),(v_2,\xi_2)]=(\{v_1,v_2\},[\xi_1,\xi_2]) 
= ( -\om(\xi_1,\xi_2),[\xi_1, \xi_2]).$$
We thus arrive at the following commutative diagram 
$$ \begin{matrix} 
& &  &  &  V^\h & \sssmapright{\id} &V^\h & \to & \0\\  
 &   &   &  & \mapdown{inc.}  & &\mapdown{inc.} & & \\  
\0 & \to & \rad(\omega,\dd_\h\omega) & \to & \hat\ham(\h,\omega) & \sssmapright{} &V_\omega & \to & \0\\  
 &   &  \mapdown{\id} &  & \mapdown{}  & &\mapdown{\dd_\h} & & \\  
\0 & \to & \rad(\omega,\dd_\h\omega)  & \to & \ham(\h,\omega) & \sssmapright{q} & \dd_\h V_\omega & \to & \0. 
\end{matrix} $$

\begin{lemm} \label{lemm1.2}
The subspace $V_\om$ is an $\sp(\h,\om)$-submodule of $V$ 
and the $V$-valued 2-cocycle $\om$
on $\h$ restricts to a $V_\om$-valued 2-cochain on $\sp(\h,\om)$.
\end{lemm}

\begin{proof}
Let $\et\in\sp(\h,\om)$ and $v\in V_\om$ with $i_\xi\om=\dd_\h v$ for some 
$\xi \in \ham(\h,\omega)$.
Then 
$$\dd_\h({\cal L}_\eta v)
={\cal L}_\et \dd_\h v
={\cal L}_\et i_\xi\omega 
=i_{[\eta, \xi]}\omega + i_\xi {\cal L}_\eta\om 
=i_{[\et,\xi]}\om$$ 
implies ${\cal L}_\et v\in V_\om$.
That $\om(\xi,\et)\in V_\om$ for $\xi,\et\in\sp(\h,\om)$ 
follows from \eqref{dh}.
\end{proof}

\begin{proposition}\label{subalg}
The central extension $\widehat\ham(\h,\om)$ is a Lie subalgebra
of the abelian extension 
$ V_\om\rtimes_\om\sp(\h,\om)$, 
defined by the Lie bracket 
$$[(v_1,\xi_1),(v_2,\xi_2)]=(\xi_1. v_2
-\xi_2. v_1+\om(\xi_1,\xi_2),[\xi_1,\xi_2]),$$
as well as of the central extension $V_\om\x_{-\om_0} \sp(\h,\om)$,
$\om_0$ denoting the 2-cocycle $\om$ on $\sp(\h,\om)$, considered as a cocycle 
with values in the trivial
$\sp(\h,\om)$-module $V_\om$. 
\end{proposition}

\begin{proof} Both assertions follow from 
$\xi_1. v_2=-\xi_2. v_1=-\om(\xi_1,\xi_2)$ 
for $\dd_\h v_j=i_{\xi_j}\om$, $j=1,2$.
\end{proof}

In the remainder of this section we describe several examples illustrating 
the abstract context described above. We start with an almost tautological 
example. 

\begin{example} \label{ex:1.4a} Let $q \: \hat \h \to \h$ be a central 
extension of Lie algebras with kernel~$\z$. Then the adjoint action 
of $\hat\h$ on itself factors through a representation 
$\hat\ad \: \h \to \der(\hat\h)$, defined by 
$\hat\ad(q(X))(Y) = [X,Y]$. Therefore $V := \hat\h$ carries a natural 
$\h$-module structure. Moreover, the Lie bracket on $\hat\h$ defines an 
invariant $2$-cocycle $\omega \in Z^2(\h,V)^\fh$, determined by 
$$ \omega(q(X),q(Y)) := -[X,Y] = q(Y).X. $$
Clearly, $V^\fh = \z(\hat\fh)$ is the center of the Lie algebra $\hat\h$ 
and $V_\omega = V$ follows from the fact that 
$i_{q(X)} \omega = \dd_\h X$ holds for each $X \in \hat\h$. This in 
turn implies that the ``Poisson bracket'' on $V = \hat\h$ is 
$$ \{ X,Y\} = - \omega(q(X),q(Y)) = [X,Y], $$
so that we recover the Lie bracket on $\hat\h$. This shows that 
any central Lie algebra extension can be written as some $V_\omega$, 
associated to an invariant $2$-cocycle. 

We also obtain $\fh = \sp(\h,\omega) = \ham(\h,\omega)$ and the corresponding 
central extension is 
\begin{align*}
\hat\ham(\h,\omega) 
&=\{ (v,X) \in V \times \fh \: i_X \omega = \dd_\h v \} 
=\{ (Y,X) \in \hat\h \times \fh \: i_X \omega = i_{q(Y)}\omega \} \\
&=\{ (Y,X) \in \hat\h \times \fh \: q(Y) - X \in \ker\hat\ad \} \\
&=\{ (Y,q(Y)) \: Y \in \hat\h\} + \{ (0,X) \: X \in \ker \hat\ad\} 
\cong \hat\h \oplus q(\z(\hat\h)). 
\end{align*}
If, in addition, 
$\z = \z(\hat\h)$, then $q(\z(\hat\h)) = \{0\}$, and we simply obtain 
$\hat\h \cong \hat\ham(\h,\omega)$, considered as a central extension 
of $\h = \ham(\h,\omega)$. 
\end{example}

\begin{examples} \label{ex:1.4} (a) If $(M,\omega)$ is a finite-dimensional connected 
symplectic manifold, 
$\h = {\cal V}(M)$ and $V = C^\infty(M,\R)$, then 
$\sp(\h,\omega) = \sp(M,\omega)$ is the Lie algebra of symplectic vector fields 
on $M$, $\ham(\h,\omega) = \ham(M,\omega)$ is the Lie algebra of hamiltonian vector fields 
on $M$, 
and $\rad \omega = \rad(\omega,\dd_\h\omega) = \{0\}$ 
implies $\h_\omega = \h$. In this case 
$$C^1(\h,V)_\omega = \{ i_\xi \omega \: \xi \in \h \} = \Omega^1(M,\R) $$
is the space of all smooth $1$-forms on $M$ and the Lie bracket on this space 
is given by 
$$ [\alpha_1,\alpha_2] = {\cal L}_{\xi_1} \alpha_2 - i_{\xi_2}{\cal L}_{\xi_1}\omega $$
for $\alpha_j = i_{\xi_j} \omega$, $j=1,2$. We also have 
$V^\h = \R$ (the constant functions), and $\hat\ham(\h,\omega) \cong (C^\infty(M,\R), 
\{\cdot,\cdot\})$ is a central Lie algebra extension of $\ham(M,\om)$. 

(b) If $\omega$ is a only a $2$-form on $M$, $\h = {\cal V}(M)$, and 
$$ \sp(M,\omega) := \{ X \in {\cal V}(M) \: {\cal L}_X \omega = 0, 
i_X \dd\om = 0\}. $$
Considering $V = C^\infty(M,\R)$ as above, we put 
$$ \ham(M,\omega) := \{ X \in {\cal V}(M) \: {\cal L}_X \omega = 0, 
(\exists f \in V)\, i_X \om = \dd f\}. $$
Then $\h_\omega$ may be a proper Lie subalgebra of 
${\cal V}(M)$. In this case the space of admissible functions 
$V_\omega = \{ f \in V \: 
(\exists \xi \in \ham(M,\omega))\ \dd f = i_\xi \omega\}$ 
is an associative subalgebra of $V$ because 
$\dd f = i_\xi \omega$ and $\dd g = i_\eta \omega$ imply 
$$ \dd(fg) = f \dd g + g \dd f = i_{f\eta} \omega + i_{g\xi}\omega 
= i_{f\eta + g\xi}\omega,$$ 
and 
\begin{align*}
{\cal L}_{f\eta}\omega 
&= \dd(i_{f\eta}\omega) 
= \dd(f i_\eta\omega) 
= \dd f \wedge i_\eta\omega + f \cdot \dd(i_\eta\omega)\\
&= \dd f \wedge i_\eta\omega + f \cdot \dd^2 g 
= i_\xi\omega  \wedge i_\eta\omega 
= - {\cal L}_{g\xi}\omega 
\end{align*}
show that $f\eta + g \xi \in \ham(M,\omega)$. 
We thus obtain on $V_\omega$ the structure of a commutative Poisson algebra 
by $\{f_1, f_2\} := \omega(\xi_2, \xi_1)$ for $i_{\xi_j} \omega = \dd f_j$ 
(cf.\ \cite{Gra85}). Hamiltonian actions in this context are 
studied in \cite{DTK07}. Our definition of $\sp(M,\omega)$ does not 
coincide with the definition given there, where the 
condition $i_X\dd \omega =0$ is omitted. This leads to a larger Lie 
algebra to which the restriction of $\omega$ is not necessarily 
a cocycle. 

For the case of closed $2$-forms, the Lie bracket on $V_\omega$ already 
occurs in Fuchssteiner's paper \cite[p.~1085]{Fu82}. 
\end{examples}

\begin{example} \label{ex:1.5} It is interesting to compare 
Examples~\ref{ex:1.4} 
with the situation arising from a Poisson 
manifold $(M,\Lambda)$. Here $\Lambda$ denotes the 
bivector field defining the Poisson bracket on $C^\infty(M,\R)$ by 
$\{f,g\} := \Lambda(\dd f,\dd g)$. 
Then we associate to $\alpha \in \Omega^1(M,\R)$ the vector field $X_\alpha 
= \Lambda^\sharp(\alpha)$, defined 
by $\beta(X_\alpha) = \Lambda(\alpha,\beta)$ for $\beta \in \Omega^1(M,\R)$.
The vector fields of the form $X_{\dd f}$ are called {\it hamiltonian}.
Then 
\begin{align*}
[\alpha, \beta] 
&:= {\cal L}_{X_\alpha}\beta - i_{X_\beta}\dd\alpha 
= {\cal L}_{X_\alpha}\beta - {\cal L}_{X_\beta}\alpha - \dd(\Lambda(\alpha,\beta)) \\
&= i_{X_\alpha}\dd\beta - i_{X_\beta}\dd\alpha + \dd(\Lambda(\alpha,\beta)) 
\end{align*}
defines a Lie bracket on the space $\Omega^1(M,\R)$ of $1$-forms on $M$ for which 
the map 
$$ \Lambda^\sharp \: \Omega^1(M,\R) \to {\cal V}(M), \quad \alpha \mapsto X_\alpha $$
is a homomorphism of Lie algebras (\cite{Fu82}, Thm.~1). We then have in particular 
$$ [\dd f, \dd g] = {\cal L}_{X_{\dd f}}\dd g 
= \dd(i_{X_{\dd f}}\dd g) = \dd\big(\Lambda(\dd f, \dd g)\big) = \dd\{f,g\}, $$
so that the exterior derivative 
$$ \dd \: (C^\infty(M,\R), \{\cdot,\cdot\}) \to (\Omega^1(M,\R),[\cdot,\cdot]) $$
is a homomorphism of Lie algebras. 

Using $\Lambda^\sharp$, we obtain the structure of an 
$\Omega^1(M,\R)$-module on $C^\infty(M,\R)$. 
Then the space $\X^p(M)$ of sections of the bundle $\Lambda^p(T(M))$ 
can be viewed as a space of $C^\infty(M,\R)$-valued Lie algebra 
$p$-cochains on $\Omega^1(M,\R)$. 
We thus obtain a subcomplex of the corresponding Chevalley--Eilenberg complex, 
on which the differential 
is given by $f \mapsto -[\Lambda,f]$, where 
$$[\cdot,\cdot] \: \X^p(M) \times \X^q(M) \to \X^{p+q-1}(M)$$ 
denotes the 
Schouten--Nijenhuis bracket (\cite{Vai94}, Prop.~4.3). 
In particular, the relation $[\Lambda,\Lambda] = 0$ implies that 
$\Lambda \in Z^2(\Omega^1(M,\R), C^\infty(M,\R))$ is a Lie algebra 
$2$-cocycle. 

With $\h := \Omega^1(M,\R)$, $V := C^\infty(M,\R)$ and $\omega := -\Lambda$ 
we are now in the setting described above. Then 
$\sp(\h,\omega)$ is a Lie subalgebra containing the Lie subalgebra of 
closed $1$-forms. Indeed, the claim being local, it suffice to prove it
for exact 1-forms, so that it follows from the relation 
\begin{align*}
X_{\dd f}(\Lambda(\dd g, \dd h))
&= \{f,\{g,h\}\}
= \{\{f,g\},h\} + \{g,\{f,h\}\}\\
&= \Lambda({\cal L}_{X_{\dd f}}\dd g, \dd h) + \Lambda(\dd g, 
{\cal L}_{X_{\dd f}}\dd h). 
\end{align*}
Since $\Lambda^\sharp$ is a homomorphism of Lie algebras, 
its kernel $\rad\omega$ is an ideal, so that $\h_\omega = \h$. 

For each $\alpha \in \Omega^1(M,\R)$ and $f \in V$ we have 
$$ (\dd_\h f)(\alpha) = X_\alpha.f = \dd f(X_\alpha) = \Lambda(\alpha, \dd f) 
= - (i_{\dd f}\Lambda) (\alpha) = -\alpha(\Lambda^\sharp(\dd f)), $$
showing that $V= V_\omega$, and the corresponding Lie algebra 
structure on $V$ coincides with the given Poisson bracket. We also note that 
\begin{align*}
V^\fh 
&= \{ f \in V \: (\forall \alpha \in \h)\ \dd f(X_\alpha) = 0\} 
= \{ f \in V \: (\forall \alpha \in \h)\ \alpha(X_{\dd f}) = 0\} \\
&= \{ f \in V \: X_{\dd f} = 0\} 
= \z(C^\infty(M,\R)) 
\end{align*}
is the center of the Poisson--Lie algebra $C^\infty(M,\R)$, 
so that 
$$ \dd_\h V = \Lambda^\sharp(\dd C^\infty(M,\R))
 \cong C^\infty(M,\R)/\z(C^\infty(M,\R)) $$
is the set $\ham(M,\Lambda)$ 
of hamiltonian vector fields of $(M,\Lambda)$. This leads to the 
short exact sequence 
$$ \0 \to V^\fh = \z(C^\infty(M,\R)) \to 
V = C^\infty(M,\R) \to
\dd_\h(V) = \ham(M,\Lambda) \to \0, $$
and the map $q \: \ham(\Omega^1(M,\R),\Lambda) \to \dd_\h(V)$ corresponds 
to the homomorphism 
$$ \Lambda^\sharp \: \ham(\Omega^1(M,\R),\Lambda) \to \ham(M,\Lambda), \quad 
\alpha \mapsto X_\alpha. $$
\end{example}

\begin{example} \label{ex:1.6} (cf.\ \cite{DKM90}, \cite{DV91}) In 
the context of non-commutative geometry, 
one considers the following situation: 
$A$ is an associative, possibly non-commutative, algebra and 
$\h$ is a Lie algebra acting by derivations on $A$. 
Then $(C^\bullet(\h,A),\dd_\h)$ is a differential graded algebra 
which can be considered as a variant of the exterior algebra 
$(\Omega^\bullet(M,\R),\dd)$ of a smooth manifold $M$. 
In the abstract context, symplectic forms on $A$ correspond 
to elements of $Z^2(\h,A)$. 

(a) Matrix algebras are particularly simple examples. For 
$A = M_n(\R)$, we consider $\h := \der(A) \cong \fsl_n(\R)$. 
Then $H^p(\h,A)$ vanishes for $p = 1,2$ by the Whitehead Lemmas. 
The Lie bracket 
$$ \omega(x,y) = [x,y] $$
defines an $\h$-invariant $A$-valued $2$-cocycle on $\h$ and since $H^1(\h,A)$ vanishes, we 
have 
$$ \h = \sp(\h,\omega) = \ham(\h,\omega). $$
The inclusion $\eta \: \h \to A$ satisfies $\dd_\h\eta = \omega$, and 
$i_x\omega = \ad x \: \h \to A,$
so that, for $x \in \h$,  $a_x := -x \in A$ 
satisfies $\dd_\h a_x = \ad x = i_x \omega$. 
We thus get $A_\omega = A$ and $A^\h = \R \1$. The Poisson bracket on 
$A_\omega = A$ coincides with the commutator bracket: 
$$  \{a,b\} = \omega(-a,-b) = [a,b].$$

(b) More generally, for any locally convex associative algebra 
$A$, the Lie algebra $\h := A/\z(A)$ acts by inner derivation on $A$,  
and the Lie bracket 
$\omega([a],[b]) := [a,b]$ is an element $\omega \in Z^2(\h,A)$. 
It is a coboundary 
if and only if the central Lie algebra extension 
$$ \0 \to \z(A) \to A \to A/\z(A) \to \0 $$
splits.\begin{footnote}{A typical example where this is not the case 
is the algebra $A = B(H)$ of bounded operators on an infinite-dimensional 
complex Hilbert space because its center $\C \1$ is contained in the 
commutator algebra (cf.\  Cor.\ 2 to Probl.\ 186 in \cite{Ha67}).} 
\end{footnote}
Here we have $\h = \ham(\h,\omega)$, 
$A = A_\omega$, $A^\h = \z(A)$, and 
$$ \{a,b\} = [a,b] $$
implies that $\hat\ham(\h,\omega) \cong A$ (as a Lie algebra). 

For any continuous linear splitting $\sigma \: \h \to A$ of the 
quotient map, we have 
$$\dd_\h \sigma(x,y) 
= [\sigma(x),\sigma(y)] - [\sigma(y),\sigma(x)] - \sigma([x,y]) 
= 2\omega(x,y) - \sigma([x,y]), $$
so that $\tilde\omega := \omega - \dd_\h\sigma$ is a $2$-cocycle 
equivalent to $\omega$, as an $A$-valued cocycle, 
but whose values lie in the trivial $\h$-module $\z(A)$.  

(c) Another, closely related example, arises for 
$\h := {\cal V}(M)$, $M$ a smooth manifold and 
the algebra $A := C^\infty(M,M_d(\R))$ for some $d > 0$. 
Then each closed $2$-form $\omega \in \Omega^2(M,\R)$ defines an 
$A$-valued $2$-cocycle on $\h$ because we may identify 
$C^\infty(M,\R)$ with the center of $A$. 
\end{example}

\subsection*{Examples arising from differential forms} 

\begin{example} \label{ex:1.7} 
Let $M$ be a finite-dimensional paracompact smooth manifold, 
$\z$ be a Mackey complete locally convex space and 
$$V := \oline\Omega^p(M,\z) := \Omega^p(M,\z)/\dd \Omega^{p-1}(M,\z). $$
We write $[\alpha] = \alpha + \dd \Omega^{p-1}(M,\z)$ for the elements 
of this space. 
In view of de Rham's Theorem, 
the subspace $\dd \Omega^{p-1}(M,\z)$ of $\Omega^p(M,\z)$ 
consists of all closed $p$-forms for which the integrals over all 
smooth singular $p$-cycles vanish. Therefore it is closed and thus 
$V$ inherits a natural Fr\'echet topology, turning it into a topological 
module of the Lie algebra $\h := {\cal V}(M)$, acting by 
$$ X.[\alpha] = [{\cal L}_X\alpha] 
= [i_X\dd\alpha + \dd i_X\alpha] 
= [i_X\dd\alpha]. $$

For any $(p+2)$-form $\tilde\omega \in \Omega^{p+2}(M,\z)$, 
we now obtain a Lie algebra $2$-cochain $\omega \in C^2({\cal V}(M),V)$ 
by 
$$ \omega(X,Y) := [i_Y i_X \tilde\omega]. $$
We may also define 
Lie subalgebras of {\it symplectic}, resp., {\it hamiltonian vector fields} 
on $M$ by 
$$ \fsp(M,\tilde\omega) := \{ X  \in {\cal V}(M) 
\: {\cal L}_X \tilde\omega = 0 = i_X \dd \tilde\omega\} $$
and 
$$ \ham(M,\tilde\omega) := \{ X \in \sp(M,\tilde\omega) \: 
i_X\tilde\omega \in \dd\Omega^p(M,\z)\}. $$ 

If $\tilde\omega$ is non-degenerate in the sense that 
$i_v \tilde\omega(p) \not=0$ for any non-zero $v \in T_p(M)$ 
and $\dd\tilde\omega=0$, then 
the pair $(M,\tilde\omega)$ is called a {\it multisymplectic manifold}. 
For more details on this class of manifolds and 
some of its applications, we refer to 
\cite[p.23]{GIMM04},  \cite{CIL99} and \cite{Ma88}. 

To understand the connection with our abstract algebraic 
setup, the following simple observation is quite useful: 
 
\begin{lemma} \label{lemma:new}
If a differential form $\alpha\in\Omega^{p+1}(M,\z)$ on the finite-dimensional 
manifold $M$ has the property
that $i_X\alpha$ is exact for any vector field $X\in{\cal V}(M)$,
then $\alpha=0$.
\end{lemma}

\begin{proof} The case $p = \dim M$ is trivial, so that we may assume that 
$p < \dim M$. 
Let $X\in{\cal V}(M)$. Then there exists $\theta\in\Omega^{p-1}(M,\z)$
with $\dd\theta=i_X\alpha$. For any $f\in C^\infty(M,\R)$ the $p$-form 
$f\dd\theta=i_{fX}\alpha$ is exact. In particular $\dd f\wedge\dd\theta=0$
for any smooth function $f$, hence $\dd\theta=0$. Then $\alpha$ vanishes
because the vector field $X$ was arbitrary.
\end{proof}

As an immediate consequence of the preceding lemma, we see that 
$X.[\alpha] = [i_X \dd \alpha]=0$ holds for all $X \in {\cal V}(M)$ 
if and only if $\dd\alpha=0$, so that 
$$V^\h= H^p_{\rm dR}(M,\z)$$
(cf.\ \cite[Lemma~23]{Ne06a}). 

\begin{proposition} \label{prop:new1} 
For the Lie algebra $2$-cochain $\omega \in 
C^2({\cal V}(M), \oline\Omega^p(M,\z))$, the following assertions hold: 
\begin{description}
\item[\rm(a)] $\omega$ is a $2$-cocycle if and only if $\tilde\omega$ is 
closed. 
\item[\rm(b)] $\sp(\fh,\omega) = \sp(M,\tilde\omega)$. 
\item[\rm(c)] $\ham(\fh,\omega) = \ham(M,\tilde\omega)$. 
\end{description}
\end{proposition}

\begin{proof} First we derive some useful formulas. 
For $X,Y,Z \in {\cal V}(M)$, we have 
$\dd\omega(X,Y,Z)
= (i_X\dd\omega)(Y,Z) 
= ({\cal L}_X\omega)(Y,Z) - \dd(i_X\omega)(Y,Z).$
Further 
\begin{eqnarray}\label{eq:act-mod}
({\cal L}_X\omega)(Y,Z)
&=&{\cal L}_X \omega(Y,Z) - \omega([X,Y],Z) - \omega(Y,[X,Z]) \notag\\
&=&[{\cal L}_X i_Z i_Y \tilde\omega - i_Z i_{[X,Y]}\tilde\omega - 
i_{[X,Z]}i_Y \tilde\omega]
=[i_Z i_Y {\cal L}_X \tilde\omega],      
\end{eqnarray}
and 
\begin{align*}
\dd(i_X\omega)(Y,Z) 
&= Y.\omega(X,Z) - Z.\omega(X,Y) - \omega(X,[Y,Z])\\
&= [{\cal L}_Y i_Z i_X \tilde\omega - {\cal L}_Z i_Y i_X \tilde\omega 
- i_{[Y,Z]}i_X \tilde\omega] \\
&= [i_Z {\cal L}_Y i_X \tilde\omega - i_Z \dd(i_Y i_X \tilde\omega)] 
= [i_Z i_Y \dd(i_X \tilde\omega)]. 
\end{align*}
This leads to 
\begin{equation}
  \label{eq:d-rel}
\dd\omega(X,Y,Z)
=[i_Z i_Y ({\cal L}_X \tilde\omega- \dd(i_X \tilde\omega))] 
=[i_Z i_Y i_X \dd\tilde\omega]. 
\end{equation}

{}From \eqref{eq:d-rel} and 
Lemma~\ref{lemma:new} we now immediately derive (a).
More precisely, we derive from Lemma~\ref{lemma:new} that for 
$X \in {\cal V}(M)$ the relation 
$i_X \dd\omega = 0$ is equivalent to 
$i_X \dd \tilde\omega = 0$ and from \eqref{eq:act-mod} that 
${\cal L}_X\omega=0$ is equivalent to 
${\cal L}_X \tilde\omega = 0$. This proves (b). 

To verify (c), we first note that for any $X\in\ham(M,\tilde\omega)$ and 
$\theta \in \Omega^p(M,\z)$ with $\dd\theta = i_X\tilde\omega$ and 
$Y \in {\cal V}(M)$, we have 
$$ Y.[\theta] = [{\cal L}_Y\theta] 
= [i_Y\dd\theta] = [i_Yi_X\tilde\omega] = \omega(X,Y) 
= (i_X\omega)(Y), $$
so that $X\in\ham(\h,\omega)$.
For the converse, let $X\in\ham(\h,\omega)$ and 
$[\theta] \in V$ with $\dd_\h[\theta]=i_X\omega$.
Then $i_Y(\dd\theta-i_X\tilde\omega)$ is exact for any vector field $Y$,
hence $\dd\theta = i_X\tilde\omega$ by Lemma \ref{lemma:new} and 
thus $X\in\ham(M,\tilde\omega)$.
\end{proof}

\begin{remark} The last  part of the 
preceding proof shows in particular that 
the Lie algebra of admissible vectors in $V = \oline\Omega^p(M,\z)$ 
is 
$$ V_\omega
=\{ [\theta] \in V \: (\exists X \in \ham(M,\tilde\omega))\ 
\dd \theta = i_X \tilde\omega\}$$ 
and its Lie bracket is 
$$ \{[\theta_1],[\theta_2]\} := \omega(X_2, X_1)
= [i_{X_1} i_{X_2} \tilde\omega]
 \quad \mbox{ for } \quad 
\dd\theta_j = i_{X_j}\tilde\omega. $$
If $\theta_i$ is closed, then we may take $X_i = 0$, which shows that 
$H^p_{\rm dR}(M,\z)$ is central in the Lie algebra $V_\omega$. 
In view of Proposition~\ref{prop:new1}, we thus obtain the central 
Lie algebra extension 
$$ \hat\ham(M,\tilde\omega) := 
\{([\theta], X) \in V_\omega \times \ham(M,\tilde\omega) \: 
\dd \theta = i_X \tilde\omega \} $$
of $\ham(M,\tilde\omega)$ by $H^p_{\rm dR}(M,\z)$, where the corresponding 
cocycle is the restriction of $-\omega$ to $\ham(M,\tilde\omega)$. 
\end{remark}
\end{example}


\section{The momentum map}\label{sec:2} 

As in the previous section, let 
$V$ a topological $\h$-module and $\omega \in C^2(\h,V)$. 
In addition, we now consider a continuous homomorphism of Lie algebras 
$$\ze:\g\to\ham(\h,\om). $$
Then we obtain via $\ze$ a topological $\g$-module structure on $V$ defined by 
$X. v := \zeta(X). v$ for $X \in \g, v \in V$ 
and the subspace $V_\om$ is a $\g$-submodule (Lemma~\ref{lemm1.2}).

\begin{definition}
The pullback 
$\hat\g_{\rm cen} := \{(v,X)\in V\x\g : i_{\ze(X)}\om=\dd_\h v\}$ 
by $\ze$ of the central extension $\hat\ham(\h,\omega)$ in 
\eqref{hatham} is 
a central extension
\begin{equation}\label{hatg}
\0\into V^\h\to\hat\g_{\rm cen}\onto\g\to \0.
\end{equation}
It can also be viewed more 
directly as the pullback of the central extension 
$V_\omega$ of $\dd_\h(V_\omega)$ by the homomorphism 
$$ \zeta_V \: \g \to \dd_\h(V), \quad X \mapsto i_{\zeta(X)}\omega. $$
\end{definition}

A continuous linear map 
\begin{equation}\label{j}
J:\g\to V_\om \quad \mbox{with} \quad
\dd_\h (J(X))=i_{\ze(X)}\om \quad \mbox{ for } \quad  X\in\g 
\end{equation}
is called a {\it momentum map} for $\ze$.

Momentum maps are in one-to-one correspondence with 
continuous linear sections 
of this extension because any continuous linear section 
$s \: \g \to \hat\g_{\rm cen}$ is of the form 
$s(X) = (J(X),X)$ for some momentum map $J$ and vice versa. 
For any such section, we obtain a $2$-cocycle by 
\begin{align*}
\ta_J(X,Y):= \ta(X,Y)&:=[s(X),s(Y)]-s([X,Y])\\
&=(\{J(X),J(Y)\},[X,Y])-(J([X,Y]),[X,Y])\\
&=( X. J(Y)-J([X,Y]),0)\in V^\h, 
\end{align*}
satisfying 
$$ \hat\g_{\rm cen} \cong V^\h \times_\tau \g. $$

\begin{lemma} \label{lem:3}  
For a momentum map $J \: \g \to V_\omega$, the following are 
equivalent: 
\begin{description}
\item[\rm(i)] $J$ is a Lie algebra homomorphism. 
\item[\rm(ii)] $J$ is $\g$-equivariant.
\item[\rm(iii)] $\tau_J = 0$. 
\item[\rm(iv)] $s = (J,\id_\g) \: \g \to \hat\g_{\rm cen}$ 
is a homomorphism of Lie algebras. 
\end{description}
\end{lemma}

\begin{proof} (i) $\Leftrightarrow$ (ii): 
The $\g$-equivariance of the momentum is equivalent to 
$X. J(Y)=J([X,Y])$ for $X,Y\in\g$. 
Hence the assertion follows from 
\begin{equation}
  \label{eq:mombrack}
 \{J(X),J(Y)\}=\om(\ze(Y),\ze(X))=(\dd_\h J(Y))(\ze(X))=X. J(Y). 
\end{equation}

\nin (ii) $\Leftrightarrow$ (iii) $\Leftrightarrow$ (iv) 
follow from the definition and the above formula for $\tau$. 
\end{proof} 

\begin{example} In the situation of Example~\ref{ex:1.4a}, 
where $V = \hat\fh$ was a central extension of $\fh$, we also have 
$\fh = \ham(\fh,\omega)$. Then an equivariant momentum map 
$J \: \fh \to \hat\fh$ is the same as the splitting of the 
central Lie algebra extension $\hat\fh$ of $\fh$ by $\fz$. 
\end{example}

For two choices $J,J'\in C^1(\g,V_\om)$ of momentum maps for $\ze$
the difference $J-J'$ has values in $V^\h$, and this leads directly to 
the following: 

\begin{proposition} 
The cohomology class $[\ta_J]\in H^2(\g, V^\h)$ does not depend on the choice
of the momentum map~$J$. It is the obstruction to the existence of 
a $\g$-equivariant momentum map for 
$\ze:\g\to\ham(\h,\om)$. In particular, the central extension 
$\hat\g_{\rm cen}$ splits if and only if an equivariant momentum map exists. 
\end{proposition}

\begin{remark}\label{jhat} 
If we replace $\g$ by the central extension $\hat\g_{\rm cen}$ 
and $\zeta$ by the 
homomorphism $\hat\ze:\hat\g_{\rm cen}\to\ham(\h,\om), 
\hat\ze(v,X) = \zeta(X)$, we obtain the $\g$-equivariant momentum map 
$$ \hat J: \hat\g_{\rm cen}= V^\h\x_\ta\g \to V, \quad (v,X) \mapsto J(X)+v.$$
Indeed, since $v \in V^\h$, we have 
$\dd_\h(\hat J(v,X)) = \dd_\h J(X) = i_{\zeta(X)}\omega$, 
and the equivariance of $\hat J$ follows from 
\begin{align*}
X.\hat J(v,Y)&-\hat J(\hat\ad(X)(v,Y))
=X.(J(Y)+v)-\hat J(\ta(X,Y),[X,Y])\\
&=X.J(Y)- J([X,Y])-\ta(X,Y))=0, 
\end{align*}
for $v\in V^\h$ and $X,Y\in\g$.
\end{remark}

The pullback cochain $\om_\g:=\ze^*\om\in C^2(\g,V)$ is a $\g$-invariant 2-cocycle with values in the $\g$-module $V_\om$
because $\ze(\g) \subeq \sp(\h,\om)$.
With 
\eqref{eq:mombrack} 
it can be expressed by the momentum map as $\om_\g(X,Y)=Y.J(X)$.

\begin{lemm} \label{lem:cobound} In $C^2(\g,V_\om)$ we have 
$\ta_J=\dd_\g J+\om_\g$. 
\end{lemm} 

\begin{proof} This follows from the relation 
  \begin{align*}
(\dd_\g J)(X,Y)&=X. J(Y)-Y. J(X)-J([X,Y])
=\ta_J(X,Y)-\om_\g(X,Y)
  \end{align*}
for $X,Y\in\g$.
\end{proof}

\begin{remark} \label{rem:newrem1} 
If $\om=\dd_\h\al$ is a coboundary and $J$ a momentum map, the 
linear map $f_J=J+\ze^*\al\in C^1(\g,V)$ has the property that for all
$X\in\g$, 
\begin{align*}
{\cal L}_{\ze(X)}\al 
&= i_{\ze(X)}\dd_\h \alpha + \dd_\h(i_{\zeta(X)}\alpha)
= i_{\ze(X)}\omega + \dd_\h(\zeta^*\alpha(X)) \\
&= \dd_\h(J(X) + \zeta^*\alpha(X))
= \dd_\h(f_J(X)),  
\end{align*}
and Lemma~\ref{lem:cobound} immediately yields $\ta_J=\dd_\g f_J$.
In general, for any $f\in C^1(\g,V)$ with ${\cal L}_{\ze(X)}\al
=\dd_\h\big(f(X)\big)$ for all $X\in\g$, we have 
$(f - f_J)(\g) \subeq V^\fh$, so that 
$[\dd_\g f] = [\dd_\g f_J] = [\tau_J]$ in $H^2(\g,V^\h)$. 
\end{remark}

\begin{definition}
Let $G$ be a connected Lie group with Lie algebra $\g$ 
and $\rh_V$ a linear action of $G$ on the $\h$-module $V$.
The $G$-action $\rh_V$ is called an {\it abstract hamiltonian action}
for the continuous 2-cochain $\om\in C^2(\h,V)$
if the derived $\g$-action on $V$ factors through $\ham(\h,\om)$.
This means there exists a Lie algebra homomorphism $\ze:\g\to\ham(\h,\om)$
with $X. v=\ze(X). v$.
When the Lie algebra homomorphism $\ze$ is given,
we call $\rh_V$ an {\it abstract hamiltonian $G$-action for~$\ze$}. 
\end{definition}

One sees immediately that $V^G=V^\g \supeq V^\h$. 
The $\g$-invariant pullback cocycle $\om_\g\in Z^2(\g,V_\om)$
is also $G$-invariant:
\begin{equation}\label{invari}
g.\om_\g(X,Y)=\om_\g(\Ad(g)X,\Ad(g)Y).
\end{equation}

Considering the natural $G$-action on $C^1(\g,V)$ by 
$(g. c)(Y)=g. c(\Ad(g)^{-1}Y)$
and its infinitesimal version 
$(X. c)(Y)=X. c(Y)-c([X,Y]),$
the Lie algebra 2-cocycle $\ta_J\in Z^2(\g,V^\h)$ satisfies 
$\ta_J=\dd_\g J\in Z^1(\g,C^1(\g,V^\h))$. As a $C^1(\g,V)$-valued cocycle 
it is a coboundary, since $J\in C^1(\g,V)$, but in general not as a 
$C^1(\g,V^\h)$-valued 1-cocycle.

We define a $C^1(\g,V)$-valued group cocycle by 
$$\ka=\dd_G J:G\to C^1(\g,V), \quad \kappa(g)=g. J-J.$$ 

\begin{lemma} \label{lem:kappa-val} 
$\ka(g)(X)\in V^\h$ for all $g\in G$ and $X\in\g$.
\end{lemma}

\begin{proof} We have already seen above that $G$ acts trivially on 
$V^\g \supeq V^\h$. To see that 
all the maps $\kappa(g)$ have values in $V^\h$, pick 
$\xi \in \h$. It suffices to show that the function 
$$F \: G \to V, \quad g \mapsto \xi . \kappa(g)X = \xi.\big((g. J-J)(X)\big) $$
is constant because it vanishes in $\1$. 
Since $G$ is connected, it suffices to see that for each 
$g \in G$ and $Y \in \g$ we have 
$$ 0 = T_g(F)(g.Y) = \xi .\big((g.(Y.J))(X)\big). $$
Since $(Y.J)(Z) = \tau(Y,Z) \in V^\h$ for $Y,Z \in \g$, and 
$(g.(Y.J))(X) = \break g . (Y.J)(\Ad(g)^{-1}X)$, this follows 
from the triviality of the action of $G$ on~$V^\fh$. 
\end{proof}

The lemma implies that $\ka$ is a $C^1(\g,V^\h)$-valued 1-cocycle on $G$.
It measures the failure of the momentum map $J$ to be $G$-equivariant because 
$$
g.\big(\ka(g^{-1})(X)\big)=J(\Ad(g)X)-g. J(X).
$$

\begin{proposition}\label{kappa}
The object measuring the failure of the momentum map $ J\in C^1(\g,V_\om)$ 
to be $G$-equivariant is the group 1-cocycle 
\begin{align*}
\ka:G\to C^1(\g, V^\h),\quad\ka(g)=g. J - J.
\end{align*}
Its cohomology class $[\ka]\in H^1(G,C^1(\g, V^\h))$ does not depend
on the choice of $J$. It is the obstruction for the existence of 
a $G$-equivariant momentum map for $\ze$.
\end{proposition}

\begin{proof}
The cohomology class $[\ka]\in H^1(G,C^1(\g,V^\h))$ does not 
depend on the choice of the momentum map.
Indeed, for two momentum maps $ J$ and $ J'$, we have 
$J-J'\in C^1(\g,V^\h)$ and
the corresponding group 1-cocycles $\ka$ and $\ka'$ satisfy 
$\kappa - \kappa' = \dd_G(J - J')$.

If $\ka$ is a $C^1(\g, V^\h)$-valued 1-coboundary on $G$,
then there is an element $c\in C^1(\g, V^\h)$ such that 
$\ka(g)=g. c-c$. Then $ J-c$ is a $G$-equivariant momentum map because, 
by definition, $\ka(g)=g. J-J$.
\end{proof}

\begin{proposition} \label{prop:2.7} 
If $\rho_V$ is an abstract Hamiltonian $G$-action for $\zeta$, 
then the 
adjoint action of $\g$ on $\hat\g_{\rm cen} = V^\h \times_\tau \g$ 
integrates to 
a smooth $G$-action, given by 
$$\hat\Ad(g)(v,X):= \big(v+ \kappa(g)(\Ad(g)X), \Ad(g)X\big) 
=\big(v-\ka(g^{-1})(X),\Ad(g)X\big), $$
for $(v,X)\in\hat\g_{\rm cen}$ and $g\in G$. With respect to this 
action, the momentum map 
$$ \hat J \: \hat\g_{\rm cen} \to V, \quad 
(v,X) \mapsto J(X)+v $$
is $G$-equivariant. 
\end{proposition}

\begin{proof} First we recall from Lemma~\ref{lem:kappa-val} that 
$\kappa(g)(\Ad(g)X) \in V^\fh$, so that $\hat\Ad(g)$ defines a continuous 
linear automorphism of $\hat\g_{\rm cen}$. 
{}From the $G$-invariance of $\omega$ and the relation 
$$ \tau(X,Y) = X.J(Y) - J([X,Y]) = \omega_\g(Y,X) - J([X,Y]), $$
we derive 
$(g.\tau-\tau)(X,Y) 
=  (J-g.J)([X,Y])
= \dd_\g(\kappa(g))(X,Y),$
so that 
$$ \hat\Ad(g)(v,X) = \big(v + \kappa(g)(\Ad(g)X), \Ad(g)X\big) $$
is a Lie algebra automorphism (\cite[Lemma~II.5]{Ne06c}). 
Further, the cocycle property of $\kappa$ implies that 
$\hat\Ad$ defines a smooth action of $G$ on $\hat\g_{\rm cen}$, 
and the corresponding derived action is 
$$ \hat\ad(Y).(v,X) = ((Y.J)(X), [Y,X]) 
= (\tau_J(Y,X), [Y,X]) = [(0,Y), (v,X)], $$
which is the $\g$-action on $\hat\g_{\rm cen}$ induced 
by the adjoint action of $\hat\g_{\rm cen}$. 
Since any representation of the connected group $G$ is determined 
by its derived representation 
(cf.\ \cite[Rem.~II.3.7]{Ne06b}), 
$\hat\Ad$ is the unique smooth action of $G$ on $\hat\g_{\rm cen}$, 
integrating the adjoint action $\hat\ad$ of $\g$ on $\hat\g_{\rm cen}$. 
Since $\kappa$ has values in $V^G$, we obtain the relation 
$$ \kappa(g)\circ \Ad(g) = g^{-1}.\kappa(g) 
= J - g^{-1}.J = - \kappa(g^{-1}). $$

Finally, the $G$-equivariance of $\hat J$ with respect to the 
$G$-action follows from the connectedness of $G$ and the equivariance 
with respect to the $\g$-action (Remark~\ref{jhat}). 
\end{proof}

We also note that the description of $\hat\g_{\rm cen}$ as 
$V^\fh \times_\tau \g$, yields an identification of the affine 
space 
$$ {\cal A} := \{ f \in \Hom(\hat\g_{\rm cen},V^\fh) \: f\res_{V^\fh} = \id_{V^\fh}\} $$
with its translation space $C^1(\g,V^\fh)$, and the $G$-action on 
${\cal A}$ induced by $\hat\Ad$ thus corresponds to the affine action 
on $C^1(\g,V^\fh)$, defined by 
\begin{equation}
  \label{eq:affine}
 (g * \alpha)(X) := \alpha(\Ad(g)^{-1}X) - \kappa(g)(X). 
\end{equation}

\begin{example} (cf.\ \cite{Ko70}) 
Let us take a closer look at the prototypical example for our 
setup. Let $(M,\omega)$ be a connected presymplectic manifold, i.e., 
$\omega$ is a closed $2$-form,  and 
$G \to \Ham(M,\omega)$ a hamiltonian action of the connected Lie group $G$ 
on $M$, where the infinitesimal action is denoted 
$\ze:\g\to\ham(M,\omega)$.

For $\h = {\cal V}(M)$ and $V= C^\infty(M,\R)$, we then have 
$V^\h = \R$, and a momentum map $J \: \g \to V$, with 
$i_{\ze(X)}\om=\dd(J(X))$ corresponds to a map 
$\mu \: M \to \g^*$ with $\mu(m)(X) = J(X)(m)$, i.e., 
$\mu$ is a momentum map for the hamiltonian $G$-action on $(M,\om)$.

The function $\ka(g)(X)=(g. J-J)(X)$ on $M$ is constant and the so 
obtained map $\ka:G\to\g^*$
is a group 1-cocycle whose cohomology class is the obstruction for the 
existence of a $G$-equivariant momentum map $\mu:M\to\g^*$. 
In any case there is an affine $G$-action on $\g^*$, defined by 
$a_g(\al)=\Ad^*(g)\al-\ka(g)$ for which $\mu$ is $G$-equivariant.
The $\g$-equivariance of the momentum map 
$\hat J \: \hat\g_{\rm cen} \to V$ 
implies its $G$-equivariance, so that the corresponding 
map $\hat\mu:M\to\hat\g_{\rm cen}^*$ is $G$-equivariant. 
By identifying $\g^*$ with the affine subspace 
$\{1\} \times \g^* \subeq (\hat\g_{\rm cen})^*$, we obtain the 
affine action~\eqref{eq:affine}: $g * \alpha := \Ad^*(g)\alpha - \kappa(g)$. 
\end{example}

\subsection*{An equivariant cohomology picture of the momentum map} 

Let $G$ be a connected Lie group with Lie algebra $\g$. 
A {\it $G^\star$-module} is a topological super vector space $\Omega$, endowed 
with a smooth $G$-action $\rho \: G \to \Aut(\Omega)$ 
by automorphisms, a continuous odd derivation $\dd_\Omega$ commuting with 
the $G$-action, and a $G$-equivariant map 
$\iota \: \g \to \End_1(\Omega)$
such that $\iota$, together with the derived representation  
${\cal L}_X := d\rho(X)$ 
turns $(\Omega,\dd_\Omega)$ into a differential graded $\g$-module. 

If $(\Om,\dd_\Omega)$ is a $G^\star$-module, then 
the corresponding {\it Cartan complex} is the space 
$$\Pol(\g,\Omega)^G := \bigoplus_{n \in \N_0} \Pol^n(\g,\Omega)^G $$
of continuous $G$-equivariant polynomial maps $\g \to \Omega$, endowed 
with the differential
$$ \dd_G(f)(X) := \dd_\Omega(f(X)) + \iota_X f(X) \quad \mbox{ for } X \in \g, 
f \in \Pol(\g,\Omega)^G. $$
If $\Omega$ is $\Z$-graded, then the natural 
grading on the space of polynomials is defined by 
$$\Pol(\g,\Omega)_d^G = \bigoplus_{2n+k=d} \Pol^n(\g,\Omega^k)^G. $$
Note that $\dd_G^2(f)(X) = {\cal L}_X(f(X))= 0$ vanishes because of the 
invariance, so that we obtain a cochain complex 
(cf.\ \cite{GS99}, Sect.~4.2). 
Its cohomology is called the $G$-equivariant cohomology
of the $G^\star$-module $\Om$.

Elements of degree $2$ in $\Pol(\g,\Omega)^G$ have the form 
$$ f = \omega + J, \quad \omega \in (\Omega^2)^G = \Pol^0(\g,\Omega^2)^G, 
\quad J \in \Lin(\g,\Omega^0)^G = \Pol^1(\g,\Omega^0)^G, $$
and the condition $\dd_G f = 0$ is equivalent to 
\begin{equation}
  \label{eq:equimo}
\dd_\Omega \omega = 0 \quad \mbox{ and } \quad 
\iota_X \omega = \dd_\Omega J(X)\quad \mbox{ for}  \quad X\in \g. 
\end{equation}

The same remains true for $G$-equivariant momentum maps
for Lie algebra 2-cocycles in the setup developed above. 
If $\h$ and $V$ are both smooth $G$-modules, so that the 
$\h$-action on $V$ is $G$-equivariant  
and the corresponding $\g$-actions come from a continuous 
homomorphism $\zeta \: \g \to \fh$, then 
the Cartan-Eilenberg complex $(C^\bullet(\h, V),\dd_\h)$
is a $G^\star$-module. Further, a $2$-cochain for the corresponding 
Cartan complex is of the form $\om+J$ for
$\om\in C^2(\h, V)^G$ and $J\in C^1(\g,V)^G$. By \eqref{eq:equimo}, the 
relation $\dd_G(\omega + J) =0$ is equivalent to $\omega$ being a 
$2$-cocycle and $J$ an equivariant momentum map for~$\omega$. 



\section{Central Lie group extensions} \label{sec:3}

Let $V$ a topological $\h$-module, $\omega \in C^2(\h,V)$ and 
$\ze:\g\to\ham(\h,\om)$ a continuous homomorphism of Lie algebras. 
Any momentum map $J$ for $\ze$ provides us with a
continuous $V^\h$-valued 2-cocycle 
$\ta(X,Y)=X. J(Y)-J([X,Y])$
on $\g$, defining the central extension $\hat\g_{\rm cen} = V^\fh \times_\tau 
\g$ from \eqref{hatg}. 
The pullback $\om_\g=\ze^*\om$ 
\begin{equation}\label{omg}
\om_\g(X,Y):=\omega(\zeta(X), \zeta(Y)) =\ze(Y). J(X) = Y. J(X)
\end{equation} 
is a $V_\om$-valued 2-cocycle on $\g$
and defines an abelian extension $\hat\g_{\rm ab} := 
V_\om\rtimes_{\om_\g}\g$ of $\g$ by $V_\om$.
In view of  Proposition \ref{subalg}, 
$\hat\g_{\rm cen}$ is a Lie subalgebra of $\hat\g_{\rm ab}$. 

Now we assume, in addition, that $V$ is Mackey complete 
and $\rho_V \: G \to \GL(V)$ is an abstract hamiltonian action for $\ze$. 
Then all period integrals are defined (cf.\ Appendix~\ref{sec:6}), and we 
assume that the period group $\Pi_{\om_\g}$ is discrete, so that 
the Lie algebra $\hat\g_{\rm ab}$ integrates to an abelian Lie group 
extension, and we may compare the groups corresponding to the Lie algebras 
$\hat\g_{\rm cen}$ and $\hat\g_{\rm ab}$.

\begin{proposition}\label{vanish}
  \begin{enumerate}
  \item[\rm(i)] The cocycles $\ta\in Z^2(\g,V^\h)$ and $\om_\g\in Z^2(\g,V_\om)$
have the same period group $\Pi_\ta=\Pi_{\om_\g} \subset V^\h$. 
\item[\rm(ii)] The flux homomorphisms $F_\ta  \: \pi_1(G) \to H^1(\g, V^\h)$ 
and $F_{\om_\g} :\pi_1(G)\to H^1(\g,V_\om)$ both vanish. 
  \end{enumerate}
\end{proposition}

\begin{proof} (i) In view of Lemma~\ref{lem:cobound}, 
$[\tau] = [\omega_\g]$ in 
$H^2(\g,V_\om)$, so that their period 
homomorphisms coincide by \cite[Thm.~7.2]{Ne04}. 

(ii) The relation $\tau - \omega_\g = \dd_\g J$ 
(Lemma~\ref{lem:cobound}) implies 
$\tau^{\rm eq} - \omega_\g^{\rm eq} = \dd J^{\rm eq}$ in $\Omega^2(G,V)$ 
(cf.\ Appendix~\ref{sec:6}), 
so that we also derive with 
${\cal L}_{X_r}J^{\rm eq}=d\rho_V(X)\circ J^{\rm eq}$ for each loop $\gamma$ 
in $G$ 
\begin{align*}
\int_\ga i_{X_r}\ta^{\rm eq} - \int_\ga i_{X_r}\om_\g^{\rm eq}
&=\int_\ga i_{X_r}(\dd J^{\rm eq}) 
=\int_\ga {\cal L}_{X_r}J^{\rm eq}  
= d\rho_V(X) .  \int_\ga J^{\rm eq}, 
\end{align*}
and since this is a $V$-coboundary, for the 
inclusion $\iota \: V^\h \to V_\om$, we have 
$$ \iota_* \circ F_\ta = F_{\omega_\g} \: \pi_1(G) \to H^1(\g,V_\om). $$
Hence it suffices to see that $F_\tau$ vanishes, but this follows 
from the existence of the smooth $G$-action on 
$\hat\g_{\rm cen}$, integrating the adjoint action of 
$\g$ (Proposition~\ref{prop:2.7} and \cite{Ne02}, Prop.~7.6). 
\end{proof}

The next theorem 
provides a central Lie group extension similar to \eqref{kostant} 
in the introduction. 

\begin{theorem}\label{thm:one} Suppose that 
$\rho_V$ is an abstract hamiltonian $G$-action on $V$ for
the Lie algebra homomorphism $\ze:\g\to\ham(\h,\om)$ 
such that the period group $\Pi := \Pi_{\omega_\g} 
= \Pi_\tau$ is discrete. Then the following assertions hold: 
\begin{description}
\item[\rm(1)] There exists a central Lie group extension $\hat G_{\rm cen}$ 
of $G$ by $Z := V^\h/\Pi$, integrating 
the Lie algebra $\hat \g_{\rm cen} = V^\h \times_\tau \g$ from~\eqref{hatg}. 
\item[\rm(2)] The quotient Lie group 
$\hat G_{\rm ab} := (V_\omega/\Pi \rtimes \hat G_{\rm cen})/\oline \Delta_Z$ 
by the antidiagonal 
$\oline\Delta_Z  :=\{ (z,z^{-1}) \: z \in Z \}$ 
is an abelian Lie group extension of $G$ by the 
smooth $G$-module $V_\om/\Pi$, integrating the Lie algebra 
$\hat \g_{\rm ab} := V_\om \rtimes_{\om_\g} \g$. 
\end{description}
\end{theorem}

\begin{proof} (1) The 2-cocycle $\ta\in Z^2(\g,V^\h)$ 
has discrete period group $\Pi$ and
vanishing flux homomorphism $F_\ta$ (Proposition \ref{vanish}). 
Hence the central Lie algebra extension 
$\hat\g_{\rm cen}=V^\h\x_\ta\g$ integrates 
to a central Lie group extension $\hat G_{\rm cen}$ 
of $G$ by $Z = V^\h/\Pi$ (Theorem \ref{integrate}).

(2) We know from Remark~\ref{rem:baer-ex} in Appendix~\ref{sec:6} below that 
$\hat G_{\rm ab}$ is a Lie group extension of $G$ by $V_\omega/\Pi$. 
Its Lie algebra is isomorphic to 
$(V_\om \rtimes \hat \g_{\rm cen})/\oline \Delta_{V^\fh}$, 
and this is the extension of $\g$ by $V_\omega$, defined by the 
cocycle $\tau$, which is equivalent to $\omega_\g$ 
(Lemma~\ref{lem:cobound}). 
\end{proof}

\begin{example}(Differential $(p+2)$-forms) \label{ex:3.3} 
Let $\rho \: G \to \Ham(M,\tilde\omega)$ 
be a smooth hamiltonian action in the context of Example~\ref{ex:1.7}.
Then the induced $G$-action on $V=\oline\Omega^p(M,\z)$
is an abstract hamiltonian action for the $V$-valued 2-cocycle $\om$ on ${\cal V}(M)$.
Let $\Gamma_{\tilde\omega} = \int_{H_{p+2}(M)} \tilde\omega \subeq \z$ be 
the group of periods of $\tilde\omega$. Then 
the Period Formula in \cite[Thm.~3.18]{Ne08} shows that 
the image of the period map 
$$ \per_{\omega_\g} \: \pi_2(G) \to H^p_{\rm dR}(M,\z) = V^\fh $$
is contained in 
$$ \Big\{ [\theta] \: \int_{H_p(M)} \theta \subeq \Gamma_{\tilde\omega} \Big\} 
\cong \Hom(H_p(M), \Gamma_{\tilde\omega}), $$
and this group is discrete whenever $\Gamma_{\tilde\omega}$ is discrete 
and $H_p(M)$ is finitely generated (which is the case if $M$ is compact). 
\end{example}

In the following section we  take a closer look at the case 
$p = 2$, where $\tilde\omega$ is a closed $2$-form on $M$. 

\begin{example} (Continuous inverse algebras) 
If, in the setting of \break Example~\ref{ex:1.6}(b), 
 $A$ is a Mackey complete continuous inverse algebra, i.e., 
its unit group $A^\times$ is open and the inversion is continuous, 
then $A^\times$ is a locally exponential Lie group, 
$Z(A)^\times = Z(A^\times)$ is a Lie subgroup, and 
$G := A^\times/Z(A)^\times$ also carries a Lie group structure 
(cf.\ \cite{GN08} for all that). The action of $G$ 
on $A$ by inner automorphisms 
is hamiltonian. Now the existence of the Lie group extension 
$$ \1 \to Z(A)^\times \to A^\times \to G \to \1 $$
implies that the corresponding period group $\Pi_\omega$ is contained 
in $\pi_1(Z(A)) = \{ z \in Z(A) \: \exp z = \1\}$ 
(cf.\ \cite{Ne02}, Prop.~V.11), hence is in particular discrete. 
\end{example}

\begin{example} (Poisson manifolds) We consider the situation 
arising for a Poisson manifold $(M,\Lambda)$ (Example~\ref{ex:1.5}). 
Suppose that $\rho \: G \to \GL(V)$ is an abstract hamiltonian 
action, corresponding to the Lie algebra homomorphism 
$\zeta \: \g \to \ham(\Omega^1(M,\R),\Lambda)$. 
Then $\g$ acts on $V = C^\infty(M,\R)$ by 
$$\zeta_V \: \g \to \ham(M,\Lambda), X \mapsto i_{\zeta(X)}\Lambda 
= \Lambda^\sharp(\zeta(X)) $$
which 
implies that $\rho(G) \subeq \Ham(M,\Lambda) \subeq \Aut(M,\Lambda)$, 
so that we obtain a Hamiltonian action of $G$ on $(M,\Lambda)$. 
The central extension $\hat\g_{\rm cen}$ is the pullback of
the central extension of $\ham(M,\Lambda) \cong \dd_\h(V)$ 
given by the Poisson--Lie algebra $V = C^\infty(M,\R)$: 
$$ \hat\g_{\rm cen} \cong \{(f,X) \in C^\infty(M,\R) \times \g \: 
\zeta_V(X) = - \Lambda^\sharp(\dd f)\}, $$
a central extension of $\g$ by $\z(C^\infty(M,\R)) = V^\fh$. 

We are interested in criteria for the corresponding period group 
$\Pi_\tau \subeq \z(C^\infty(M,\R))$ to be discrete. 
Let us assume that $(M,\Lambda)$ admits a quantization line bundle 
$q \: \bL \to M$, i.e., $\bL$ is a complex line bundle on which 
we have a covariant derivative $\nabla$ for which the operators 
$$ \hat f.s := \nabla_{X_{\dd f}} s + 2 \pi i f \cdot s $$
define a homomorphism of Lie algebras 
$$ C^\infty(M,\R) \to \End(\Gamma\bL), \quad f \mapsto \hat f. $$
A characterization of Poisson manifolds for which 
such bundles exist is given in \cite{Vai91} (see also \cite{Vai94} 
and \cite{Hu90}). It is equivalent to the existence of a closed 
$2$-form $\lambda$ representing an integral cohomology class 
and a vector field $A$ for which 
the bivector field $\lambda^\sharp$ defined by 
$\lambda^\sharp(\alpha,\beta) := \lambda(\alpha^\sharp,\beta^\sharp)$ 
satisfies 
$$ \Lambda + L_A \Lambda = \lambda^\sharp. $$
Then $\lambda$ is the curvature of the pair $(\bL,\nabla)$, 
hence a $2$-cocycle describing the Lie algebra extension 
$$ \0 \to \End_{C^\infty(M,\R)}(\Gamma\bL) \cong C^\infty(M,\R) 
\to \dend(\Gamma\bL) \to {\cal V}(M) \to \0, $$
where $\dend(\Gamma \bL)$ is the Lie algebra of 
derivative endomorphisms, i.e., those endomorphisms $D \in \End(\Gamma\bL)$ 
for which there exists a vector field $X \in {\cal V}(M)$, such that 
$[D,\nabla_X]$ is multiplication with a smooth function 
(cf.\ \cite{Ko76}). 

We conclude that the restriction of the $C^\infty(M,\R)$-valued 
cocycle $\zeta_V^*\lambda$ is equivalent to $\zeta^*\tau$, hence 
leads to the same period homomorphism 
$$ \per \: \pi_2(G) \to V^\fh \subeq V. $$
Since the existence of the line bundle $\bL$ with curvature 
$\lambda$ implies that all periods of $\lambda$ are integral, 
the discussion in Example~\ref{ex:3.3} (specialized to $p =0$), 
implies that the period group of $\zeta^*\tau$ is discrete, 
so that Theorem~\ref{thm:one} implies the existence of a 
central Lie group extension $\hat G_{\rm cen}$ with 
Lie algebra $\hat\g_{\rm cen}$ acting by bundle automorphisms 
on~$\bL$. The crucial difference to the symplectic case is that 
in this situation the Lie algebra $\hat\g_{\rm cen}$ acts on 
$\bL$ by vector fields not necessarily preserving the connection, 
resp., the covariant derivative. 
\end{example}

\begin{problem} Suppose that $(M,\Lambda)$ is a Poisson manifold 
(Example~\ref{ex:1.5}). When is the Lie algebra $\Omega^1(M,\R)$ 
integrable in the sense that it is the Lie algebra of some infinite-dimensional Lie group? Since it is an abelian extension of the 
Lie algebra $\im\Lambda^\sharp$ of vector fields, one would 
like to prove first that this Lie algebra is integrable and 
then try to use the results on abelian extensions in \cite{Ne04}, 
but the integrability of such Lie algebras of vector fields 
defined by integrable distributions is a 
difficult problem which is still open (cf.\ \cite{Ne06b}, Problem~IV.13). 

We can ask the same question for the Poisson--Lie algebra 
$C^\infty(M,\R)$, which is a central extension of the Lie algebra 
of Hamiltonian vector fields. Under which conditions does any of 
these two Lie algebras integrate to a Lie group? If $(M,\Lambda)$ 
is compact symplectic, then $\ham(M,\Lambda)$ always does 
and its central extension $C^\infty(M,\R)$ does at least if the 
cohomology class of the symplectic form has a discrete period group 
(cf.\ Example~\ref{ex:3.3}). 
\end{problem}


\section{From symplectic to hamiltonian actions} \label{sec:4}

Let $\z$ be a Mackey complete locally convex space, 
$M$ a connected smooth manifold (possibly infinite-dimensional) and 
$\omega \in \Omega^2(M,\z)$ a closed $2$-form.  
Under the rather weak assumption 
that  the group $S_\omega  := \int_{\pi_2(M)} \omega$ of {\sl spherical} periods of $\omega$ 
is discrete, we use a bypass through the simply connected 
covering $q_M \: \tilde M\to M$ to associate to a smooth symplectic 
action 
$$G \to \Sp(M,\omega) := \{ \phi \in \Diff(M) \: \phi^*\omega = \omega\}$$ 
of a connected Lie group $G$ a hamiltonian action of the simply connected 
covering group $\tilde G$ on $\tilde M$ 
and further a corresponding central Lie group extension. 

We now turn to the details. 
Let $Z = \z/\Gamma_Z$, where $\Gamma_Z \subeq \z$ is a discrete 
subgroup. We write $q_Z \: \z \to Z$ 
for the quotient map and assume that the 
group $S_\omega$ of 
spherical periods of $\omega$ is contained in $\Gamma_Z$. 
Whenever $S_\omega$  is discrete, this is in particular 
satisfied if we put $\Gamma_Z := S_\omega$. 

We write $\pi_1(M)$ for the group of deck transformations of 
$\tilde M$ over $M$, acting from 
the left. Then 
$\tilde \omega := q_M^*\omega$ is a closed $\z$-valued $2$-form 
on $\tilde M$, and since the natural homomorphism 
$\pi_2(\tilde M) \to \pi_2(M)$ is an isomorphism, it has the same group 
$S_\omega$ of spherical periods. Moreover, the Hurewicz homomorphism 
$\pi_2(\tilde M) \to H_2(\tilde M)$ is an isomorphism, so that 
all periods of $\tilde\omega$ are contained in $\Gamma_Z$. 
If, in addition, $M$ is smoothly paracompact, then this implies the 
existence of a pre-quantum 
principal $Z$-bundle $q_P \: P \to \tilde M$ with a connection 
$1$-form $\theta \in \Omega^1(P,\z)$ whose curvature is $\tilde\omega$, i.e., 
$q_P^*\tilde\omega = \dd\theta$ (cf.\ \cite{Br93}). 

The following theorem is a slight generalization of 
Kostant's Theorem concerning finite-dimensional manifolds and 
the case $Z = \T$ (\cite{Ko70}, Prop.~2.2.1). 

\begin{theorem} \label{thm:quanto-ext} We have an abelian group extension 
$$ \1 \to C^\infty(\tilde M,Z) \cong \Gau(P) 
\to \Aut(P) \to \Diff(\tilde M)_{[\tilde\omega]} \to \1 $$
and a central extension 
$$ \1 \to Z \to \Aut(P,\theta) \to \Sp(\tilde M,\tilde\omega) \to \1. $$
\end{theorem}

\begin{proof} Smooth $Z$-bundles over $\tilde M$ are classified by the 
group 
$$ H^2(\tilde M,\Gamma_Z) \cong \Hom(\pi_2(M),\Gamma_Z) \into 
H^2_{\rm dR}(\tilde M,\z) \cong \Hom(\pi_2(M),\z). $$
Therefore $\phi \in \Diff(\tilde M)$ lifts to a bundle 
isomorphism of $P$ if and only if 
$[\phi^*\tilde\omega] = [\tilde\omega]$, which in turn is equivalent to 
$\phi^*P \cong P$ as $Z$-bundles. This leads to the abelian extension, 
where $C^\infty(\tilde M,Z) \cong \Gau(P)$ acts on $P$ 
by $\phi_F(p) := p.F(q_P(p))$. 

Any quantomorphism $\tilde\phi \in \Aut(P,\theta)$ factors 
through an element of the group 
$\Sp(\tilde M,\tilde\omega)$. 
To see that, conversely, each element $\phi$ of $\Sp(\tilde M, \tilde\omega)$ 
lifts to a 
quantomorphism of $P$, we first note that the preceding paragraph 
yields the existence of a lift $\tilde\phi$ to some bundle automorphism. 
Then $\tilde\phi^*\theta - \theta$ can be written as $q_P^*\alpha$ for some 
$\alpha \in \Omega^1(\tilde M,\z)$, and we have 
$$ q_P^*\dd\alpha = \dd q_P^*\alpha 
= \tilde\phi^*(\dd\theta) - \dd\theta 
= \tilde\phi^*q_P^*\tilde\omega - q_P^*\tilde\omega 
= q_P^*(\phi^*\tilde\omega - \tilde\omega), $$
so that $\alpha$ is closed if $\phi \in \Sp(\tilde M,\omega)$. 
As $H^1_{\rm dR}(\tilde M,\z)$ vanishes, there exists a 
smooth function $f \: \tilde M \to \z$ with $\dd f =  \alpha$. 
For $F := q_Z \circ f \in C^\infty(\tilde M,Z)$ and the corresponding 
gauge transformation $\phi_F \in \Gau(P)$, we then have 
$$ \phi_F^*\theta - \theta = q_P^*\dd f = q_P^*\alpha 
= \tilde\phi^*\theta - \theta, $$
so that $\tilde\phi \circ \phi_F^{-1} \in \Aut(P,\theta)$ is a lift 
of $\phi$. Now the observation that 
$\Aut(P,\theta)$ intersects $\Gau(P)$ in $Z$ 
leads to the desired central extension. 
\end{proof}

Now let $G$ be a connected Lie group and 
$\rho \: G \to \Sp(M,\omega)$ a homomorphism defining a smooth action of 
$G$ on $M$ preserving $\omega$. Then there exists a 
unique smooth action $\tilde\rho \: \tilde G \to \Sp(\tilde M,\tilde\omega)$ 
of the universal covering group $\tilde G$ of $G$ on $\tilde M$. 
Further, $\tilde\omega = q_M^*\omega$ is invariant under 
$\pi_1(M)$. Let 
$\tilde\Sp(M,\omega) := N_{\Sp(\tilde M,\tilde\omega)}(\pi_1(M))$ 
denote the normalizer of $\pi_1(M)$ in $\Sp(\tilde M,\tilde\omega)$, 
which coincides with the set of all 
lifts of elements of $\Sp(M,\omega)$ to $\tilde M$. 
We thus obtain a short exact sequence 
$$ \1 \to \pi_1(M) \to \tilde\Sp(M,\omega) \to \Sp(M,\omega) \to \1 $$
with $\tilde\rho(\tilde G) \subeq \tilde\Sp(M,\omega)$. 

Since $\tilde G$ is connected and 
$\ham(\tilde M,\tilde\omega) = \sp(\tilde M,\tilde\omega)$, the 
action of $\tilde G$ on $\tilde M$ is hamiltonian, and we 
derive from \cite{NV03}, Prop.~1.12 and Thm.~3.4, 
that the pullback $\hat G_{\rm cen} := \tilde\rho^*\Aut(P,\theta)$ 
is a central Lie group extension of $\tilde G$ by $Z$, 
i.e., $\hat G_{\rm cen}$ carries a Lie group structure which 
is a principal $Z$-bundle over~$\tilde G$. A Lie algebra 2-cocycle
representing the corresponding central Lie algebra extension 
$\hat\g_{\rm cen}$ of $\g$ by $\z$ 
is given in terms of the derived action 
$$\tilde\zeta \: \g \to \ham(\tilde M,\tilde\omega)$$ 
of $\tilde\rho$, resp., 
the derived action $\zeta \: \g \to \sp(M,\omega)$ of $\rho$ by 
$$ \tau(X,Y) 
= -\tilde\om(\tilde\ze(X),\tilde\ze(Y))(\tilde m_0) 
= -\om(\ze(X),\ze(Y))(m_0), $$
where $\tilde m_0 \in \tilde M$ and $m_0 \in M$ are points with 
$q_M(\tilde m_0) = m_0$. 
This further leads 
to an abelian Lie group extension 
$$ \1 \to C^\infty(M,Z)_0 \to \hat G_{\rm ab} \to \tilde G \to \1, $$
integrating the cocycle $\zeta^*\omega \in Z^2(\g,C^\infty(M,\z))$. 

Since the kernel of the universal covering homomorphism 
$q_G \: \tilde G \to G$ is the abelian group $\pi_1(G)$, 
the group $\hat G_{\rm cen}$ can also be viewed as an extension 
of $G$ by the group $\hat\pi_1(G)$, which is the inverse image of 
$\ker q_G \cong \pi_1(G)$ in $\hat G_{\rm cen}$. 
As a central extension 
$$ \1 \to Z \to \hat\pi_1(G) \to \pi_1(G) \to \1 $$
of an abelian group, this group is $2$-step nilpotent. Since 
$Z$ is divisible, all extensions of $\pi_1(G)$ which are abelian groups 
are trivial, 
so that $\hat\pi_1(G)$ is characterized by its commutator map 
$$ C \: \pi_1(G)  \times \pi_1(G) \to Z, $$
(\cite{Bro82}, Thm.~6.4). This map can be calculated by 
$$ C([\alpha],[\beta]) 
= q_Z\Big(\int_{\tilde\alpha \bullet \tilde\beta} \tau^{\rm eq}\Big), $$
where $\tau^{\rm eq} \in \Omega^2(\tilde G,\z)$ is the left invariant 
$2$-form corresponding to $\tau$ and 
$$ \tilde\alpha \bullet \tilde\beta \: [0,1]^2 \to\tilde G, \quad 
(t,s) \mapsto \tilde\alpha(t)\tilde\beta(s), $$
where $\tilde\alpha, \tilde\beta \: [0,1] \to \tilde G$ 
are lifts of the smooth loops $\alpha,\beta$ in $G$, starting in~$\1$ 
(cf.\ \cite{Ne04}, Cor.~6.5). 
Since the $2$-form $\tau^{\rm eq}$ on $\tilde G$ is the pullback 
of the corresponding form on $G$, we also have 
$$ C([\alpha],[\beta]) 
= q_Z\Big(\int_{\alpha \bullet \beta} \tau^{\rm eq}\Big) 
= q_Z\Big(\int_{\T^2} (\alpha \bullet \beta)^*\tau^{\rm eq}\Big), \quad 
\alpha \bullet \beta(t,s) := \alpha(t)\beta(s). $$
The orbit map $\rho^{m_0} \: G \to M$ is 
equivariant with respect to the left action of $G$ on itself, so that 
$(\rho^{m_0})^*\omega \in \Omega^2(G,\z)$ is the left invariant 
$2$-form on $G$ whose value in $\1$ is $-\tau$. Therefore 
$\tau^{\rm eq} = - (\rho^{m_0})^*\omega$, which leads to 
$$ C([\alpha],[\beta]) 
= q_Z\Big(-\int_{(\alpha \bullet \beta).m_0} \omega\Big),$$
where $(\alpha \bullet \beta).m_0$ 
can be considered as a smooth map $\T^2 \to M$. 

\begin{remark}  We 
have introduced $P$ as a $Z$-bundle over $\tilde M$, but we can 
also interprete it as a bundle over $M$, as follows. 
Let $\hat\pi_1(M) \subeq \Aut(P,\theta)$ denote the inverse image of the 
discrete subgroup $\pi_1(M) \subeq \Sp(\tilde M,\tilde\omega)$. 
Then we have a short exact sequence 
$$ \1 \to Z \to \hat\pi_1(M) \to \pi_1(M) \to \1, $$
and this group carries a natural Lie group structure for which 
$Z$ is an open central subgroup, $\pi_1(M)$ is its group of 
connected components, and the action of this group on $P$ is smooth. 
The orbit space of this action is $P/\hat\pi_1(M) 
\cong \tilde M/\pi_1(M) \cong M$ and since the map to $M$ has smooth 
local sections, $P$ is a smooth $\hat\pi_1(M)$-principal 
bundle over $M$. The extension $\hat\pi_1(M)$ 
splits if and only if there exists a 
$Z$-prequantum bundle $(P_M,\theta_M)$ for $(M,\omega)$ 
(cf.\ \cite{Ko70}, Prop.~2.4.1). 
\end{remark} 


\section{The algebraic essence of the Noether Theorem}\label{sec:5}

Given a hamiltonian $G$-action on a presymplectic manifold $(M,\om)$ 
($\omega$ is a closed $2$-form on $M$),
a momentum map $\mu:M\to\g^*$ and a $G$-invariant function $f$ on $M$,
Noether's Theorem says that along the trajectories of the hamiltonian vector field $X_f := X_{\dd f}$, the momentum map is constant.
An important result of Marsden and Weinstein \cite{MW74}
concerns symplectic reduction in the context where 
two connected Lie groups $G_1$ and $G_2$
act in a hamiltonian way on the presymplectic manifold $(M,\om)$
with momentum maps $\mu_1,\mu_2$, such that $\mu_2$ is constant 
along the $G_1$-orbits.
They show that $\mu_1$ is constant along the $G_2$-orbits and the two actions commute. In the present section we show that the algebraic essence of this 
result can be formulated in our abstract setup of Lie algebra 
$2$-cocycles.  

Let $V$ a topological $\h$-module, $\omega \in Z^2(\h,V)$ and 
$\ze:\g\to\ham(\h,\om)$ a continuous homomorphism of Lie algebras. 
The space $V^\g$ of $\g$-invariant vectors for the topological 
$\g$-module structure on $V$
obtained via $\ze$ contains $V^\h$.
The space of admissible $\g$-invariant vectors 
$V_\om^\g=V^\g\cap V_\om$ is a Lie subalgebra of $V_\om$ containing
$V^\h$ (cf.\ Proposition~\ref{prop:1}).

\begin{proposition}
If $\xi \in \ham(\h,\omega)$ satisfies $i_\xi \omega = \dd_\h v$ for some 
$v \in V^\g$, then $\xi. J(X)=0$ holds for all $X\in\g$.
\end{proposition}

\begin{proof}
Let $v\in V^\g_\om$ and $\dd_\h v=i_\xi\om$. Then
$$\xi. J(X)=\big(\dd_\h J(X)\big)(\xi)=i_\xi i_{\ze(X)}\om
=-i_{\ze(X)}\dd_\h v=-\ze(X). v=0 $$
for all $X\in\g$.
\end{proof}

\begin{proposition}
Let $\ze_j:\g_j\to\ham(\h,\om)$, $j=1,2$,  be two Lie algebra homomorphisms
and $J_j:\g_j\to V$, $j=1,2$, corresponding 
momentum mappings such that $J_2$ takes values in 
$V^{\g_1}$. Then $J_1$ takes values in $V^{\g_2}$ and 
$[\ze_1(\g_1),\ze_2(\g_2)] \subeq \rad \omega$.
\end{proposition}

\begin{proof}
Let $X_1\in\g_1$ and $X_2\in\g_2$. Then $\ze_1(X_1). J_2(X_2)=0$.
But $\ze_2(X_2). J_1(X_1)=\om(\ze_1(X_1),\ze_2(X_2))
=-\ze_1(X_1). J_2(X_2)=0$ and we further derive 
$0=\dd_\h (\om(\ze_1(X_1),\ze_2(X_2)))=-i_{[\ze_1(X_1),\ze_2(X_2)]}\om,$
so that $[\ze_1(X_1),\ze_2(X_2)]\in\rad \omega$.
\end{proof}


\section{Appendix: Integration of abelian Lie algebra extensions}\label{sec:6}

According to the general theory developed in \cite{Ne02} and 
\cite{Ne04}, there 
are two obstructions for the integration of a Lie algebra 
cocycle $\omega \in Z^2(\g,V)$ with values in a smooth Mackey complete 
$G$-module $V$ to a Lie group extension of $G$ by a quotient group 
$A = V/\Gamma_A$, where $\Gamma_A$ is a discrete subgroup of $V$:  
the period map and the flux homomorphism. To define these homomorphisms, 
we associate to $\alpha \in C^p(\g,V)$ the left equivariant 
$V$-valued $p$-form $\alpha^{\rm eq} \in \Omega^p(G,V)$, determined by 
$$ \alpha^{\rm eq}_\1 = \alpha, \quad 
\lambda_g^*\alpha^{\rm eq} = \rho_V(g) \circ \alpha^{\rm eq}, $$
where $\rho_V \: G \to \GL(V)$ describes the $G$-module 
structure on $V$. We write $d \rho_V \: \g \to \gl(V)$ for the 
corresponding derived representation. 

The {\it period map} is the group homomorphism 
$$\per_\om:\pi_2(G)\to V^G, \quad 
\per_\om([\si])=\int_{\SS^2}\si^*\om^{\rm eq} \quad 
\mbox{ for } \quad \si\in C^\infty(\SS^2,G).$$ 
Its image $\Pi_\omega$ 
is called the {\it period group} of $\om$. 

The {\it flux homomorphism} $F_\om:\pi_1(G)\to H^1(\g,V), [\gamma] 
\mapsto [I_\gamma^\omega]$, 
assigns to each piecewise smooth 
loop $\ga$ in $G$ based at the identity,
the cohomology class of the 1-cocycle
$$I_\ga^\omega \: \g\to V, \quad I_\gamma^\omega(X) = -\int_\ga i_{X_r}\om^{\rm eq},$$
where $X_r \in {\cal V}(G)$ denotes the right invariant vector field 
with $X_r(\1) = X$. 

\begin{theo} {\rm(\cite{Ne04}, Thm.~6.7)}\label{integrate}
For a Lie algebra 2-cocycle $\om\in Z^2(\g,V)$  with discrete
period group $\Pi_\om$ and vanishing flux homomorphism,
the Lie algebra extension $\hat\g=V \rtimes_\om \g$
integrates to an abelian Lie group extension
$$
\1\to V/\Pi_\om\into\hat G\onto G\to \1.
$$ 
\end{theo}

\begin{definition} \label{def:baer-prod} Let 
$G$ be a group, $A$ an abelian group 
which is a $G$-module and $Z \subeq A^G$ a subgroup. Then we define the 
{\it Baer product} of a central extension 
$q_c \: \hat G_c \to G$ of $G$ by $Z$ and an abelian extension 
$q_a \: \hat G_a \to G$ of $G$ by $A$ by 
$$ \hat G_c \otimes \hat G_a := \hat G/\{(z,z^{-1}) \: z \in Z \}, $$
where 
$$ \hat G := \{ (g_1, g_2) \in \hat G_c \times \hat G_a \: 
q_c(g_1) = q_a(g_2)\} $$
is the fiber product of the two extensions, which is an abelian 
extension of $G$  by the product module $Z \times A$, and 
the antidiagonal 
$$\oline\Delta_Z  :=\{ (z,z^{-1}) \: z \in Z \} \subeq Z \times A 
\subeq \hat G $$
is central in $\hat G$.
\end{definition}

\begin{remark} \label{rem:baer-ex} (a) 
On the level of cocycles, the Baer product corresponds to the natural 
action map $H^2(G,Z) \times H^2(G,A) \to H^2(G,A)$, 
induced by the multiplication map $Z \times A \to A, (z,a) \mapsto za$. 

(b) Suppose, in addition, $G$, $Z$ and $A$ above are Lie groups, where 
the action of $G$ on $A$ is smooth. 
Then the Baer product of two Lie group extensions 
$\hat G_c$ and $\hat G_a$ carries a natural structure of a Lie group 
extension of $G$ by $A$. Here we use that the antidiagonal 
$\oline\Delta_Z$ is a Lie subgroup of $Z \times A$ with 
$Z \times A \cong \oline\Delta_Z \times A$, so that the factorization 
of $\oline\Delta_Z$ defines a Lie group extension. 

(c) If the extension $\hat G_a$ is trivial, 
i.e., of the form $\hat G_a = A \rtimes G$, then we have 
in the notation of the preceding definition 
$\hat G \cong A \rtimes \hat G_c,$
where $\hat G_c$ acts on $A$ through the quotient map $q_c$, and 
this leads to 
$$ \hat G_c \otimes \hat G_a \cong (A \rtimes \hat G_c)/\oline\Delta_Z, $$
which is the abelian extension of $G$ by $A$ obtained from $\hat G_c$ 
by the natural inclusion $Z \into A$. 
\end{remark}

\nin Karl-Hermann Neeb \\ 
Mathematics Department \\ 
Darmstadt University of Technology \\ 
Schlossgartenstrasse 7 \\  
65289 Darmstadt \\ 
Germany \\ 
email: neeb@mathematik.tu-darmstadt.de \\

\nin Cornelia Vizman \\ 
Mathematics Department \\ 
West University of Timisoara \\ 
Bd. V.Parvan 4 \\ 
300223 Timisoara \\ 
Romania \\ 
email: vizman@math.uvt.ro


\begin{thebibliography}{aaaaaaa} 

\bibitem[Bro82]{Bro82} Brown, K. S., ``Cohomology of Groups,'' 
Grad. Texts Math. {\bf 87}, Springer-Verlag, 1987 

\bibitem[Br93]{Br93}
Brylinski, J. L.,
``Loop Spaces, Characteristic Classes and Geometric Quantization,''
Progress in Mathematics {\bf 107}, Birkh\"auser, 1993.

\bibitem[CIL99]{CIL99} Cantrijn, F., L. A. Ibort, and M. de L\'eon, 
{\it On the geometry of multisymplectic manifolds}, J. Aust. Math. Soc. A 
{\bf 66} (1999), 303--330 

\bibitem[DTK07]{DTK07} Di Terlizzi, L., and 
J. J. Konderak, {\it Reduction Theorems for Manifolds with Degenerate 2-form}, 
J.  Lie Theory {\bf 17} (2007), 563--582 

\bibitem[DV91]{DV91} Dubois-Violette, {\it Noncommutative differential 
geometry, quantum mechanics and gauge theory}, in 
``Differential geometric methods in theoretical physics,'' 
Lecture Notes in Physics {\bf 375}, Springer-Verlag, 1991; 13--24 

\bibitem[DKM90]{DKM90}
Dubois-Violette, M., R. Kerner and J. Madore, 
{\it Noncommutative differential geometry of matrix algebras}, 
J. Math. Phys. {\bf 31:2} (1990), 316--322 

\bibitem[FF01]{FF01} Feigin, B. L. and D. B. Fuchs, {\it Cohomologies of Lie Groups and 
Lie Algebras}, in ``Lie Groups and Lie Algebras II,'' A. L. Onishchik and 
E. B. Vinberg (Eds.), Encyclop. Math. Sci. {\bf 21}, Springer-Verlag, 2001 

\bibitem[Fu82]{Fu82}
Fuchssteiner, B., 
{\it The Lie algebra structure of degenerate Hamiltonian and bi-Hamiltonian 
systems}, Progr. Theoret. Phys. {\bf 68:4} (1982), 1082--1104. 

\bibitem[GIMM04]{GIMM04} Gotay, M. J., J. Isenberg, J. E. Marsden, and
R. Montgomery, {\it Moment Maps and Classical Fields. Part I: Covariant Field Theory}, arXiv:physics/9801019v2, August 2004 

\bibitem[GN08]{GN08} 
Gl\"ockner, H., and K.-H. Neeb, 
``Infinite Dimensional Lie Groups, Vol. I, 
Basic Theory and Main Examples,'' 
book in preparation 

\bibitem[Gra85]{Gra85} 
Grabowski, J.,  
{\it The Lie structure of $C^*$- and Poisson algebras}, 
Studia Math. {\bf 81} (1985), 259--270. 

\bibitem[GS99]{GS99} 
Guillemin, V., and S. Sternberg,
``Supersymmetry and Equivariant de Rham Theory,''
Springer Verlag, Berlin, 1999.

\bibitem[Ha67]{Ha67} Halmos, P. R., 
``A Hilbert Space Problem Book,'' Graduate Texts
in Math.\ {\bf 19}, Springer-Verlag, 1967 

\bibitem[Hu90]{Hu90} Huebschmann, J., 
{\it Poisson cohomology and quantization}, J. reine angew. Math. 
{\bf 408} (1990), 57--113 


\bibitem[Ko76]{Ko76}  Kosmann, Y., 
{\it On Lie transformation groups and the covariance of 
differential operators}, in ``Differential geometry and relativity,'' 
pp.~75--89, Math. Phys. and Appl. Math. {\bf 3}, Reidel, Dordrecht, 1976 

\bibitem[Ko70]{Ko70} 
Kostant, B.,  
{\it Quantization and unitary representations}, in 
``Lectures in Modern Analysis and Applications III,'' 
Springer Lecture Notes Math. {\bf 170} (1970),
87--208 

\bibitem[MW74]{MW74} 
Marsden, J., and A. Weinstein, 
{\it Reduction of symplectic manifolds with symmetry}, 
Rep. Mathematical Phys. {\bf 5:1} (1974), 121--130.

\bibitem[Ma88]{Ma88} 
Martin, G., {\it A Darboux theorem for multi-symplectic manifolds}, 
Lett. Math. Phys. {\bf 16} (1988), 133--138 



\bibitem[Ne02]{Ne02} Neeb, K.-H., 
{\it Central extensions of infinite-dimensional
Lie groups}, Annales de l'Inst. Fourier {\bf 52} (2002), 1365--1442 

\bibitem[Ne04]{Ne04} ---, {\it Abelian 
extensions of infinite-dimensional Lie groups}, 
Travaux Math. {\bf 15} (2004), 69--194.

\bibitem[Ne06a]{Ne06a} ---,
{\it Lie algebra extensions and higher order cocycles},
J. Geom. Symm. Phys. {\bf 5} (2006), 1--27.

\bibitem[Ne06b]{Ne06b} ---, {\em Towards a Lie theory of locally convex 
groups}, Jap. J. Math. 3rd series {\bf 1:2} (2006), 291--468. 

\bibitem[Ne06c]{Ne06c} ---, {\it Nonabelian 
extensions of topological Lie algebras}, 
Communications in Alg. {\bf 34} (2006), 991--1041 

\bibitem[Ne08]{Ne08} ---, {\it Lie groups of bundle automorphisms 
and their extensions}, to appear in 
``Trends and Developments in Infinite Dimensional Lie Theory,'' 
K.-H. Neeb and A. Pianzola, Eds., Progress in Mathematics, 
Birkh\"auser Verlag, 2008 

\bibitem[NV03]{NV03}
Neeb, K.-H. and C. Vizman,
{\it Flux homomorphisms and principal bundles over infinite dimensional
manifolds}, Monatsh. Math. {\bf 139} (2003), 309--333.

\bibitem[Vai91]{Vai91}
Vaisman, I.,  {\it On the geometric quantization of Poisson manifolds}, 
J. Math. Phys. {\bf 32:12} (1991), 3339--3345 

\bibitem[Vai94]{Vai94}
Vaisman, I.,  
``Lectures on the Geometry of Poisson Manifolds,'' 
Progr. Math. {\bf 118}, Birkh\"auser, 1994. 



\end{thebibliography}
\end{document}